\documentclass[12pt,leqno]{article}
\usepackage{amsmath,amssymb,latexsym}
\usepackage{theorem}
\allowdisplaybreaks
\def\nek{,\ldots,}
\def\CC{\mathbb C}
\def\HH{\mathbb H}
\def\NN{\mathbb N}
\def\RR{\mathbb R}
\def\TT{\mathbb T}
\def\ZZ{\mathbb Z}
\def\alp{\alpha}
\def\oalp{\overline\alpha}
\def\bet{\beta}
\def\del{\delta}
\def\Del{\Delta}
\def\eps{\epsilon}
\def\ome{\omega}
\def\intl{\int\limits}
\def\supl{\sup\limits}
\def\suml{\sum\limits}
\def\gam{\gamma}
\def\Gam{\Gamma}
\def\Sig{\Sigma}
\def\lam{\lambda}
\def\part{\partial}
\def\kap{\kappa}
\def\Ome{\Omega}
\def\MOD{{\rm mod}}
\def\const{{\rm const}}

\def\G{\mathcal G}
\def\calA{\mathcal A}

\def\calL{\mathcal L}
\def\done{{1\hskip-2.5pt{\rm l}}}

\def\tilome{{\widetilde\ome}}

\def\of{{\overline f}}
\def\oh{{\overline h}}
\def\tilf{{\widetilde f}}
\def\tilF{{\widetilde F}}
\def\tilg{{\widetilde g}}
\def\tilh{{\widetilde h}}
\def\tilM{{\widetilde M}}
\def\tilx{{\widetilde x}}
\def\tily{{\widetilde y}}

\def\tilrho{{\widetilde\rho}}

\newcommand{\pr}{\medskip\noindent\textit{Proof}. }
\newcommand{\Diff}{{\text{Diff}}}
\newcommand{\Symp}{{\text{Symp}}}
\newcommand{\Ham}{{\text{Ham}}}
\newcommand{\Flux}{{\text{Flux}}}
\newcommand{\width}{{\text{width}}}
\newcommand{\spectrum}{{\text{spectrum}}}
\newcommand{\Length}{{\text{length}}}
\newcommand{\Area}{{\text{area}}}

\def\bks{{\backslash}}
                
\def\nek{,\ldots,}

\def\tilF{{\widetilde F}}

\newtheorem{theorem}{Theorem}[subsection]
\newtheorem{lemma}[theorem]{Lemma}
\newtheorem{corollary}[theorem]{Corollary}
\newtheorem{proposition}[theorem]{Proposition}

{\theorembodyfont{\rmfamily}
\newtheorem{example}[theorem]{Example}}
{\theorembodyfont{\rmfamily}
\newtheorem{remark}[theorem]{Remark}}
{\theorembodyfont{\rmfamily}
\newtheorem{definition}[theorem]{Definition}}

\def\theequation{\thesection.\arabic{equation}}

\begin{document}
\title{Growth of maps, distortion in groups and
symplectic geometry}
\author{Leonid Polterovich
\thanks{Supported by the ISRAEL SCIENCE FOUNDATION founded
by the Israel Academy of Sciences and Humanities.}
\\
School of Mathematical Sciences\\
Tel Aviv University\\
Ramat Aviv, Israel  69978\\ }

\date{December 5 , 2001}
\maketitle


\begin{abstract}
In the present paper we study two sequences of real numbers
associated to a symplectic diffeomorphism:
\begin{itemize}
\item{}The uniform norm of the differential of its $n$-th iteration;
\item{}The word length of its $n$-th iteration, where we assume that
our diffeomorphism lies in a finitely generated group of
symplectic diffeomorphisms.
\end{itemize}
We find lower bounds for the growth rates of these
sequences in a number of situations. These bounds depend
on the symplectic geometry of the manifold rather than on the
specific choice of a diffeomorphism. They are obtained by using
recent results of Schwarz on Floer homology.
As an application, we prove non-existence of certain non-linear
symplectic representations for finitely generated groups.
\end{abstract}

\hfill\eject

\tableofcontents

\hfill\eject

\baselineskip=18pt
\section{Introduction and main results}\label{sec:intro}

\subsection{Growth and distortion}\label{subsec:growth}
Given a diffeomorphism $f$  of a smooth compact
connected manifold $M$, define its {\it growth
sequence\/}
$$\Gam_n(f)=\max \Big(\max_{x\in M}| d_xf^n|,\
\max_{x\in M}| d_xf^{-n}|\Big)\ ,\ n\in\NN\ .$$
Here $| d_xf|$  stands for the operator norm of the
differential $d_x f$  calculated with respect to a
Riemannian metric on $M$.  Though the explicit value of
$\Gam_n(f)$  depends on the choice of  metric, the
appropriately defined {\it growth type\/} is an
invariant of $f$  under conjugations in $\Diff (M)$.
Here is the definition.  Given two positive sequences
$a_n$  and $b_n$, we write $a_n\succeq b_n$ if there
exists $c>0$  so that $a_n\ge cb_n$  for all $n\in\NN$,
and $a_n\sim b_n$  if $a_n\succeq b_n$  and $b_n\succeq
a_n$. With this language the growth type of a
diffeomorphism $f$  is simply the equivalence class of
the sequence $\Gam_n(f)$.

The interest to the growth type is caused by a number
of reasons.  Using this notion, one can imitate the
fundamental trichotomy hyperbolic-parabolic-elliptic in
the context of diffeomorphisms (cf. \cite{HK}).
We say that $f$  is
hyperbolic if $\Gam_n(f)$ is growing exponentially fast,
$f$ is elliptic if $\Gam_n(f)$  is bounded and $f$ is
parabolic otherwise.  This definition of course agrees
with the classical one for M\"obius transformations
acting on the circle.  Sometimes the type of a
diffeomorphism reflects its important dynamical
features.  Here are several examples:
Existence of an invariant measure
with a positive Lyapunov exponent yields hyperbolicity.
If $f$  can be included into an action of a compact
group then it is elliptic.  Integrable systems of
classical mechanics often give rise to parabolic
diffeomorphisms with $\Gam_n(f)\sim n$  (see an example
in  1.4.C below).

An interesting problem, which goes back to D'Ambra and
Gromov \cite{DAG} is to study restrictions on the growth
type for various classes of diffeomorphisms.  In the
present paper we give some answers in the symplectic
category.  In particular, we prove a lower
bound for the growth type of symplectic diffeomorphisms
of closed symplectic manifolds with vanishing $\pi_2$
(see 1.4 below).
If $M$ is a closed oriented surface endowed with an
area form $\ome$, the lower bounds for the growth type
look as follows.  Denote by $\Symp_0(M,\ome)$  the group of
all area-preserving diffeomorphisms of $M$  isotopic to
the identity map $\done$.

\begin{theorem}\label{theo1.1A} Let $M$ be a closed
oriented surface of genus $\ge 2$.  Then
$\Gam_n(f)\succeq n$ for every $f\in
\Symp_0(M,\ome)\bks\{\done\}$.
\end{theorem}

\begin{theorem}\label{theo1.1B} {\bf \cite{PS}}.  Let
$M=\TT^2$ be the $2$-torus, and let \\ $f\in
\Symp_0(\TT^2,\ome)\bks\{\done\}$ be a symplectic
diffeomorphism with a fixed point.  Then
$\Gam_n(f)\succeq n$.
\end{theorem}

\medskip
\noindent
We refer to 1.4 and 1.5 below for further discussion
and generalizations to higher dimensions.

\medskip
From a different, more geometric, viewpoint the function
$$\log \Gam_1: \Diff (M)\to [0; +\infty)$$  is a
pseudo-norm on $\Diff (M)$.  Then the growth type of $f$
reflects the distortion of the cyclic subgroup $\{f^n\}
\subset \Diff(M)$  with respect to this pseudo-norm.
This observataion
serves as a motivation for the study
of distortion in finitely generated groups of symplectic
diffeomorphisms (see 1.6 and \S 4 below). Here is a sample result
for surfaces. We write $||f||$ for the word length of an
element of a finitely generated group.

\begin{theorem}\label{theo1.1.C}
Let $M$ be a closed oriented surface of genus $\geq 2$,
and let $\G \subset \Symp_0 (M,\ome)$ be a finitely
generated subgroup. Then $||f^n|| \succeq \sqrt{n}$
for every $f \in \G \setminus \{\done\}$.
\end{theorem}

\medskip
\noindent

Theorem \ref{theo1.1.C} and
related results
(see 1.6 and \S 4 below) have a number of applications to
the Zimmer program of studying non-linear representations
of discrete groups \cite{Zimmer}. We refer to papers
\cite{Ghys},\cite{BurMon},
\cite{FM},\cite{KM},\cite{FF} for recent exciting developments
in this direction. Our first applications deal with the group
$G$ of all (not neccesarily isotopic
to identity) smooth symplectic diffeomorphisms of a closed oriented
surface $M$ of genus $\geq 2$.
The next corollary was explained to us
by Marc Burger.

\begin{corollary}
Let $\G$ be a irreducible non-uniform lattice in
a semisimple real Lie group of real rank at least two.
Assume that the Lie group is connected, without compact
factors and with finite center.
Then every homomorphism $\G \to G$  has finite image.
\end{corollary}

\noindent
The proof is given in 1.6 below.
A prototype example of such a lattice is $SL(n,\ZZ)\subset SL(n,\RR)$
for $n \geq 3$.

Next, consider
the Baumslag-Solitar group
$$BS(q,p)=\langle a,b\mid
a^q=ba^pb^{-1}\rangle\ ,\ \text{where}\
q,p\in\ZZ\ ,\ q\not= 0\ ,\ p\not= 0\ ,\, |p| < |q|.$$

\begin{theorem}\label{cor1.1.D}
For every homomorphism $\phi: BS(q,p) \to G$ the
element $\phi(a)$ is of finite order.
\end{theorem}

\noindent
The proof is given in 4.7 below. In 1.6 the reader will
find generalizations to actions in higher dimensions.
We refer to \cite{FF} for the study of actions of $BS(q,p)$
on 1-dimensional manifolds.

This circle of problems is quite sensitive to the class
of smoothness
of diffeomorphisms in question, see for instance \cite{FS}
for results on real-analytic actions of lattices on surfaces.
Throughout
the paper we work with $C^{\infty}$-diffeomorphisms, however our results
and proofs should be valid in
the $C^1$-case (see 4.6 below for an outline of such
an extension).

Our approach to  growth and distortion is
based on some properties of the action
spectrum of Hamiltonian diffeomorphisms which were
obtained by Schwarz \cite{Schwarz} with the use of Floer
homology.
Interestingly enough, the bounds we get in higher
dimensions substantially depend on the fundamental group
of $M$.  In order to state our results we need the
notion of the symplectic filling function (cf.
\cite{Gromov},\cite{Sikorav}) which will be introduced
right now.

\subsection{Symplectic filling function}
\label{subsec:symplectic}
Let $(M,\ome)$  be a closed symplectic manifold with
$\pi_2(M)=0$.  Denote by $\tilome$ the lift of the
symplectic structure $\ome$  to the universal cover
$\tilM$ of $M$.  The condition $\pi_2(M)=0$ guarantees that
$\tilome$ is exact on $\tilM$.  Let $\calL$ be the space
of all 1-forms on $\tilM$  whose differential equals
$\tilome$. Fix any Riemannian metric, say $\rho$, on $M$
and write $\tilrho$ for its lift to $\tilM$.  Pick a
point $x\in\tilM$  and denote by $B(s)$  the Riemannian
ball of radius $s > 0$  centered at $x$.
Put
$$u(s) =\inf\limits_{\alp\in\calL}\ \supl_{z\in
B(s)}\ |\alp_z|_\tilrho\ .$$
Clearly the function $s\mapsto su(s)$  is strictly
increasing.  Let $v:(0; +\infty)\to (0; +\infty)$  be
its inverse.  We call $v$  {\it the symplectic filling
function\/} of $(M,\ome)$.  It is easy to check
(see 3.1 below) that if
$v'$ is the symplectic filling function associated to
another Riemannian metric $\rho'$  on $M$  and a base
point $x'\in\tilM$ then $c^{-1}v\le v'\le cv$  for some
$c>0$  (which is denoted by $v\sim v'$).

\begin{example}\label{exam1.2A}
Consider the standard symplectic torus
$\TT^{2n}=\RR^{2n}/\ZZ^{2n}$ with the symplectic
form $ \ome =  \sum^n_{j=1} dp_j\wedge
dq_j.$ We claim that $u(s)\sim s$, and therefore
$v(s)\sim\sqrt{s}$.  Indeed, take $x=0\in\RR^{2n}$, and
let $\rho$  be the Euclidean metric.  Then $\tilome
=d\Big(\sum^n_{j=1}p_jdq_j\Big)$, so
$$u(s)\le \sup\limits_{|p|^2+|q|^2\le s^2}\
{\Big |}\suml^n_{j=1}p_jdq_j{ \Big|}=\supl_{|p|^2+|q|^2\le
s^2}|p|=s\ .$$
On the other hand, consider the 2-disc $D(s)$ of radius
$s$  in the $(p_1,q_1)$-plane.  Note that for every
primitive $\alp\in\calL$
$${\rm length}\part D(s) \cdot \sup\limits_{z\in
B(s)}|\alp_z|\ge\intl_{\part D(s)}\alp
=\intl_{D(s)}\tilome =\pi s^2\  ,$$
and hence $$\supl_{z\in B(s)}|\alp_z|\ge\frac{\pi
s^2}{2\pi s}=\frac{s}{2}\;.$$  Thus
$u(s) \geq s/2$ and
the claim follows.
\end{example}

\begin{example}\label{exam1.2B}
Let $(M,\ome)$  be a closed oriented surface of genus
$\ge 2$.  We claim that $u(s)$  is bounded, and
therefore $v(s)\sim s$.  To prove the claim, represent
$M$  as $\HH/K$, where $\HH=\{ p+iq\in\CC\mid q >0\}$
is the hyperbolic upper half-plane, and $K$  is a
discrete group of isometries.  The hyperbolic metric
$\tilrho$  on $\HH$  is given by
$(dp^2+dq^2)/q^2$.  Assume without loss of
generality that the lift $\tilome$  of the symplectic
form coincides with the hyperbolic area form:  $\tilome
=(dp\wedge dq)/q^2$.  Note that $\tilome =d\alp$
for $\alp =  dp/q$.  Take $z=p+iq\in\HH$  and
calculate $$
|\alp_z|_{\tilrho}=\sup\limits_{(\xi,\eta):
\xi^2+\eta^2=q^2} \frac{|\eta|}{q} =1.$$
Hence $u(s)\le 1$  and the claim follows.  This example
motivates the next definition.
\end{example}

\begin{definition}\label{def1.2C}
A closed symplectic manifold $(M,\ome)$  with
$\pi_2(M)=0$  is called {\it symplectically
hyperbolic\/} if the function $u(s)$  is bounded (and
therefore $v(s)\sim s$).  For instance surfaces of genus
$\ge 2$,  their products, and, more generally, K\"ahler
hyperbolic manifolds \cite{GJDG}, are symplectically
hyperbolic.
\end{definition}

\subsection{Fixed points of symplectic diffeomorphisms}
\label{subsec:fixed}

Let $(M,\ome)$ be a closed symplectic manifold.  Denote
by $\Symp_0(M,\ome)$  the identity component of the group
of all symplectic diffeomorphisms of $M$.

\medskip\noindent
{\bf Definition.} Let $x$  be a fixed point of a
symplectic diffeomorphism \\ $f\in \Symp_0(M,\ome)$.  We say
that $x$  is {\it of contractible type\/} if there exists
a path $\{ f_t\}_{t\in [0;1]}$ of symplectic
diffeomorphisms with $f_0 = \done, f_1=f$ such that the loop
$\{ f_tx\}_{t\in [0;1]}$  is contractible in $M$.

\begin{example}\label{exam1.3A}
Every Hamiltonian diffeomorphism of a closed symplectic
manifold with $\pi_2=0$  has a fixed point of
contractible type.  This is an immediate consequence of
Floer's famous proof of the Arnold conjecture (see
\cite{Fl}; we refer to \cite{Pbook} for background
on Hamiltonian diffeomorphisms and to 2.2 below
for the definition.)
\end{example}

\begin{example}\label{exam1.3B}
Consider the standard symplectic torus
$\TT^{2n}=\RR^{2n}/\ZZ^{2n}$ endowed with the
symplectic form $dp\wedge dq$.  Let
 $f\in
\Symp_0(\TT^{2n})$ be a symplectic diffeomorphism.
We claim that {\it every} fixed point
$x$ of $f$ is of contractible
type.  Indeed, take a path $\{ f_t\}_{t\in [0;1]}$
of symplectic diffeomorphisms of $\TT^{2n}$  such that
$f_0=\done$  and $f_1=f$.  Let
$\tilf_t:\RR^{2n}\to\RR^{2n}$  be its lift to the
universal cover.  Pick a lift
$\tilx$  of $x$.  Then $\tilf_1\tilx =\tilx +e$ for some
$e\in\ZZ^{2n}$.
Consider a symplectic flow $g_t:\TT^{2n}\to\TT^{2n}$
defined by
$$g_tz=z-te\ (\MOD \; 1)\ ,\ z\in\TT^{2n}\ .$$
Put $h_t=g_tf_t$ for $t\in [0;1]$.  Since $g_0=g_1=\done$
we have $h_0=\done$ and $h_1=f$.  Further, the lift
$\{\tilh_t\}$ of $\{h_t\}$ to $\RR^{2n}$  satisfies
$\tilh_t$ $\tilx=\tilf_t\tilx -te$, and in particular
$\tilh_1\tilx =\tilx$. Hence the loop $\{h_tx\}_{t\in
[0;1]}$ is contractible on $\TT^{2n}$.  This completes
the proof of the claim.
\end{example}

\begin{example}\label{exam1.3C}
Let $M$ be a closed surface of genus $\ge 2$.  We claim
that every $f\in \Symp_0(M)$  has a fixed point of
contractible type.  Indeed, write $M=\HH/K$  were
$\HH$  is the hyperbolic upper half-plane and $K$ is
a discrete group of M\"obius transformations.  Take a
symplectic isotopy $\{f_t\}_{t\in[0;1]}$ with
$f_0=\done$, $f_1=f$  and lift it to $\HH$. We get an
isotopy $\{\tilf_t\}$.  Put $\tilf =\tilf_1$.  Clearly
it suffices to show that $\tilf$  has a fixed point on
$\HH$.  Assume on the contrary that $\tilf x\not= x$
for all $x\in\HH$.  Define a smooth vector field $\eta$
on $\HH$  as follows:  $\eta(x)$ is the unit tangent
vector to the geodesic ray $[x; \tilf(x))$. Since
$\tilf$ commutes with elements of $K$, the field
$\eta$ is $K$ - invariant and hence descends to a
unit vector field on $M$.  Since the Euler
characteristic of $M$ does not vanish we get a
contradiction. The claim follows.
\end{example}

\subsection{A lower bound for the growth type}\label{subsec:main}
We are ready now to state our main result.  Let
$(M,\ome)$  be a closed connected symplectic manifold with
$\pi_2(M)=0$.  Let $v$  be its symplectic filling
function.

\begin{theorem}\label{theo1.4A}
Let
$f\in
\Symp_0(M,\ome)\bks \{\done\}$ be a
symplectic diffeomorphism
with a fixed point of
contractible type. Then  $\Gam_n(f)\succeq v(n)$.
\end{theorem}

\medskip
\noindent
The proof is
given in \S 3 below.  Several remarks are in order.

\medskip
\noindent
{\bf 1.4.B.} As a consequence we get that $\Gam_n(f)\succeq\sqrt{n}$
if $(M,\ome)$  is the standard symplectic torus
(see \ref{exam1.2A}). Theorem \ref{theo1.1B} refines this
estimate for $n=2$.  Further, $\Gam_n(f)\succeq n$  if
$(M,\ome)$  is symplectically hyperbolic (see
\ref{exam1.2B},\ref{def1.2C}).  In particular, Theorem
\ref{theo1.1A} is an immediate consequence of
\ref{theo1.4A},\ref{exam1.3C} and \ref{exam1.2B}.
Patrice LeCalvez informed us that he can prove Theorems
\ref{theo1.1A} and \ref{theo1.1B} by a different method.
In higher dimensions, however, I am not aware of any alternative
to the symplecto-topological approach.

\medskip
\noindent
{\bf 1.4.C.} On every compact symplectic manifold $(M,\ome)$  one can
find a Hamiltonian diffeomorphism $f\not=\done$ such
that $\Gam_n(f)\sim n$.
Indeed, suppose that $\dim M=2m$.  Put
$N_\eps =\TT^m\times D^m(\eps)$, $\eps >0$  where
$\TT^m=\RR^m(q_1\nek q_m)/{\ZZ^m}$ and $D^m(\eps)=\{
p\in\RR^m:|p|\le\eps\}$.  Endow $N_\eps$  with the
standard symplectic form $\Ome =\sum^m_{i=1}dp_i\wedge
dq_i$.  It is well known that for $\eps >0$  small
enough there exists a symplectic embedding
$j:(N_\eps,\Ome)\to (M,\ome)$.  Fix such $\eps$  and
$j$, and take any function $H:D^m(\eps)\to\RR$ which
vanishes near $\part D^m(\eps)$, and is not identically
zero.  Consider the Hamiltonian flow of $H=H(p)$  on
$N_\eps$:
$$h_t(p,q)=(p,q+t\frac{\part H}{\part p}(p))\ $$
(this flow represents the simplest {\it integrable}
system of Classical Mechanics).
Obviously, $\Gam_n(h_1)=\Gam_1(h_n)\sim n$.  Now define
a Hamiltonian diffeomorphism $f:M\to M$  as follows:
$$f\equiv\done\ \text{on}\ M\bks j(N_\eps),\ \text{and}\
f=jh_1j^{-1}\ \text{on}\ j(N_\eps)\ .$$
Clearly, $\Gam_n(f)\sim n$  as required.

We conclude that the inequality $\Gam_n(f)\succeq n$ on
$\TT^2$  and on any symplectically hyperbolic manifold
is sharp and cannot be improved.

\medskip
\noindent
{\bf Open problem.}  Find a closed symplectic
manifold $M$ with $\pi_2(M)=0$ and a Hamiltonian
diffeomorphism $f\not=\done$ of $M$ which violates the
inequality $\Gam_n(f)\succeq n$.

\medskip
\noindent
{\bf 1.4.D.} The condition $\pi_2(M)=0$  cannot be removed.  Indeed
consider the standard $S^1$-action $\{f_t\}$ on $S^2$  (rotation
around the vertical axis).  For each $t$ the diffeomorphism $f_t$
is Hamiltonian, has fixed points of
contractible type at the poles, and obviously the
sequence $\Gam_n(f_t)$  is bounded.  (Still, it sounds
likely that $\Gam_n(f)\succeq n$  for every non-identical
area-preserving map $f:S^2\to S^2$  with at least 3 fixed
points.)

\medskip
\noindent
{\bf 1.4.E.} The fixed point condition of Theorem \ref{theo1.4A}
cannot be removed.  Indeed, $\{\Gam_n(f)\}$ is bounded
if $f$  is a translation of $\TT^2$.  More
sophisticated counterexamples are given by the next
theorem.

\setcounter{theorem}{5}
\begin{theorem}\label{theo1.4}
For every $\bet\in (0;1)$  there exists a
$C^\infty$--function $\psi :S^1\to\RR$ and an irrational
number $\alp$ so that the map $$f:\TT^2\to\TT^2,
(x,y)\mapsto (x+\alp,y+\psi (x))$$  satisfies the
following:
\begin{itemize}
\item [{\rm (i)}] $\Gam_n(f)\preceq n^\bet\log n$;
\item[{\rm (ii)}] $\Gam_{n_i}(f)\succeq n^\bet_i$  for a
subsequence $n_i\to +\infty$.
\end{itemize}
\end{theorem}

\medskip
\noindent
The proof is quite technical and will be given elsewhere.
Let us emphasize that the growth bound 1.4.A is
in general not true if $f$ has fixed points but none of them
is of contractible type (see Appendix).

\medskip
\noindent
{\bf 1.4.G.} Let $G$ be a group of diffeomorphisms acting on a
compact manifold $M$.  We say that an increasing
function $w:(0;+\infty)\to (0; +\infty)$, $w(s)\to
+\infty$  as $s\to +\infty$ is a {\it growth bound\/}
for $G$  if $\Gam_n(f)\succeq w(n)$  for all $f\in
G\bks\{\done\}$.

\subsection*{Examples:}
{\bf 1.4.G(i).} Suppose that $G$  is the group of Hamiltonian
diffeomorphisms of a closed symplectic manifold $M$.
Then the symplectic filling function $v(s)$ gives a
growth bound for $G$ in view of Theorem \ref{theo1.4A}.
One can take $w(s)=s$  for a symplectically hyperbolic
manifold.

\medskip
\noindent
{\bf 1.4.G(ii).} Fix a point $x_0$ on a closed symplectic manifold $M$.
Let $G$ be the identity component of the group $\{f\in
\Symp_0(M)\mid f(x_0)=x_0\}$.  The conclusions of the
previous example are still valid for $G$  in view of
\ref{theo1.4A}.

\medskip
\noindent
{\bf 1.4.G(iii).} Fix a point $x_0$ on any
compact manifold $M$.  Let $G$
be the identity component of the group $\{ f\in
\Diff_0(M)\mid f(x_0)=x_0\}$.  It is proved in \cite{PSo}
that $G$  admits no growth bound.  In fact, for every
increasing function $w:(0; +\infty)\to (0; +\infty)$
with $w(s)\to +\infty$  as $s\to +\infty$  there
exists a diffeomorphism $f\in G \bks \{\done\}$  and a subsequence
$n_i\to +\infty$ so that $\Gam_{n_i}(f)\preceq w(n_i)$.

\medskip
\noindent
{\bf 1.4.G(iv).} Let $M$ be a closed
manifold endowed with a volume form,
and let $x_0$  be a point on $M$.  Define $G$ as the
identity component of the group of all volume-preserving
diffeomorphisms $f$  with $f(x_0)=x_0$.

\medskip
\noindent
{\bf Open problem.} Does there exist a closed
manifold $M$  of dimension $\ge 3$  such that $G$
admits a growth bound?

\subsection{Growth and propagation}\label{subsec:propag}
The growth type of a diffeomorphism $f$  of a compact
manifold $M$  is related to the dynamics of its lift $\tilf$
to the universal cover $\tilM$. Fix such a lift and
consider a fundamental domain $D$ of $\tilM$.
Take a Riemannian metric $\rho$  on $M$  and write
$\tilrho$  for its lift to $\tilM$.  For purposes of our
discussion assume that $\tilf x=x$  for some $x\in D$.
Consider the quantity
$$d_n(\tilf)=\sup\limits_{z\in D}
\text{distance}_{\tilrho}(x,\tilf^nz)\ ,$$
which measures the rate of propagation of the
trajectories of $\tilf$  on the universal cover
(cf. \cite{P2},\cite{BPS}).
Obviously,
\renewcommand{\theequation}{\thesubsection.\Alph{equation}}
\begin{equation}\label{eq1.5A}
\Gam_n(f)\succeq d_n(\tilf)\ .
\end{equation}
In some situations this inequality combined with
an information about fixed points of $f$  gives rise to a
lower bound for the growth type of $f$.  Assume for
instance that there exists a point $x'\in D$
such that $\tilf x'=Tx'$, where $T$  is the deck
transformation corresponding to an element
$\alp\in\pi_1(M)$.  Then obviously
$d_n(\tilf)\succeq\|\alp^n\|$ where $\|\quad \|$  stands
for the word length in $\pi_1(M)$.  Thus
$\Gam_n(f)\succeq \|\alp^n\|$.

Let $(M,\omega)$ be a closed symplectic manifold
with $\pi_2 = 0$. Consider a Hamiltonian diffeomorphism
$f$ of $M$.
In his famous work \cite{Fl} on Arnold's conjecture
Floer proved that every Hamiltonian path $\{f_t\}$
with $f_0 = \done$ and $f_1 =f$  has a contractible $1$-periodic
orbit. Consider
the lift $\tilf$ of $f$ to the
universal cover $\tilM$  associated to the
path $\{f_t\}$.
It follows that the lift $\tilf$ does not depend
on the specific choice of a Hamiltonian path
joining the identity with $f$.
We will call $\tilf$  {\it the canonical
lift} of $f$. Note also that
a contractible $1$-periodic
orbit of $\{f_t\}$ corresponds to a fixed point of
$\tilf$.
Thus the sequence
$\{ d_n(\tilf)\}$ is well defined and its growth
type  is an invariant of the Hamiltonian diffeomorphism $f$.

We say that $f$  {\it does not propagate\/} if
$d_n(\tilf)$ is bounded.
One may have
the impression
that this property is not too useful since the estimate
\eqref{eq1.5A} becomes trivial.  Paradoxically, in
the Hamiltonian category the lack of propagation guarantees
at least linear growth of the differential.

\setcounter{theorem}{1}
\begin{theorem}\label{theo1.5B}
Let $f\not=\done$  be a Hamiltonian
diffeomorphism of a closed symplectic manifold
with $\pi_2(M)=0$.  Assume that $f$ does not propagate.
Then $\Gam_n(f)\succeq n$.
\end{theorem}

\medskip
\noindent
The proof is given in \S 3 below.

\subsection{Symplectic filling function and
distortion in finitely generated groups}
\label{subsec:filling}
Let $(M,\ome)$  be a closed connected symplectic manifold
with $\pi_2(M)=0$, and let $\G$
 be a finitely generated
subgroup
of
$\Symp_0(M,\ome)$.
Fix a system of generators in $\G$
and write $\| f\|$ for the word length of an
element $f\in\G$.  We are interested in the
{\it distortion\/} of the cyclic subgroup $\{
f^n\}\subset\G$, that is in the growth type of
the sequence $\| f^n\|$ (see Ch.~3 in \cite{Gromov}
for discussion on the distortion). Our main results
in this direction are given in the next two theorems.

\begin{theorem}\label{theo1.6A}
Let $\G\subset \Ham (M,\ome)$ be a finitely
generated subgroup. Then
$\|f^n\|\succeq v(n)$  for all $f\in \G\bks\{\done\}$.
\end{theorem}

\begin{theorem}\label{theo1.6B}
Assume that the fundamental group $\pi_1(M)$
has trivial center.
Let $\G \subset \Symp_0 (M,\ome)$ be a finitely
generated subgroup.
Then $$\| f^n\|\succeq\min
(v(n),\sqrt{n})$$  for all $f\in\G\bks
\{\done\}$.
\end{theorem}

\medskip
\noindent
The assumption on $\pi_1$ in Theorem 1.6.B cannot be removed.
For instance, the group
$\Symp_0(\TT^{2m},dp\wedge dq)$
contains a translation $f$ of a finite order,
thus $\{\| f^n\|\}$  is a bounded sequence.
Theorem \ref{theo1.1.C} is an immediate
consequence of 1.6.B.
{\rm Theorems \ref{theo1.6A}} and {\rm \ref{theo1.6B}} are
proved in \S 4.
Some remarks are in order.

\medskip
\noindent
{\bf 1.6.C.} Let $(M,\ome)$  be a symplectically hyperbolic
manifold.  Then $v(n)\sim n$ and we conclude
that every element $f\not=\done$  of a finitely
generated subgroup $\G\subset \Ham (M,\ome)$
is {\it undistorted\/}:  $\| f^n\|\sim n$. This
follows from \ref{theo1.6A} and the obvious
upper bound $\| f^n\|\preceq n$.

\medskip
\noindent
{\bf 1.6.D.}
Theorem 1.6.A can be rephrased as follows.  Let
$\G$ be an abstract finitely generated
discrete subgroup.  Assume that $g\in\G$ is
an element such that $\| g^{n_i}\|/v(n_i) \to
0$  as $n_i\to +\infty$.  Then $\phi(g)=\done$
for every homomorphism $\phi:\G\to \Ham (M,\ome)$.
Theorem 1.6.B, of course, admits a similar reformulation.

Let us illustrate this statement.
Following \cite{LMR} we call an element $x \in \G$  a {\bf U}-element if
it is of infinite order and
$$\lim \inf \frac{\log||x^n||}{\log n} = 0.$$
{\bf U}-elements appear in a number
of interesting groups. We present two examples.

\medskip
\noindent
{\bf Example 1.6.E.} Consider the Baumslag-Solitar group
$$BS(q,p)=\langle a,b\mid
a^q=ba^pb^{-1}\rangle\ ,\ \text{where}\
q,p\in\ZZ\ ,\ q\not= 0\ ,\ p\not= 0\ ,\; |p| < |q| .$$
It is known that the element
$a$  has logarithmic distortion:  $\|
a^n\|\preceq\log (n+1)$.  Here is a simple
argument which we learned from Zlil Sela.
Assume for simplicity that $q>p>0$.  Define a
function $\varphi :\NN\to\ZZ$  as follows:
$\varphi(k)$ is the integer lying in
$\Big(\frac{p}{q}k-1;\frac{p}{q}k\Big]$.  Note
that $a^{qk} =ba^{pk}b^{-1}
=ba^{q\varphi(k)+i_k}b^{-1}$, where $i_k\in[0;q-1]$.
Put $u_k=\| a^{qk}\|$.  Then $u_k\le
q+1+u_{\varphi(k)}$ for all $k\in\NN$  which
readily yields $u_k\le\const\cdot\log(k+1)$  for
all $k\in\NN$.  But $\| a^{qk+j}\|\le q-1+u_k$
for all $k\in\NN, j\in\{1;\cdots ; q-1\}$.  This
yields $\| a^n\|\preceq\log (n+1)$  as required.

\medskip
\noindent
{\bf Example 1.6.F.} Let $\G$ be an irreducible
lattice in a semisimple real Lie group of real rank at least two.
We assume that the Lie group is connected, without compact
factors and with finite center.
If $\G$ is non-uniform, a classical
result due to Kazhdan and Margulis
implies that it contains a unipotent
element of infinite order.
Lubotzky-Mozes-Raghunathan \cite{LMR} proved that
it must be a {\bf U}-element.

\medskip
\noindent
\setcounter{theorem}{6}
\begin{corollary}\label{corol1.6E}
Let $G$  be one of the following groups:
\begin{itemize}
\item[(i)] $\Ham(M,\ome)$ where $(M,\ome)$ is either
the standard symplectic $2m$-dimensional torus
or a
symplectically
hyperbolic manifold;
\item[(ii)] $\Symp_0(M,\ome)$  where $M$  is a
product of surfaces of genus $\ge 2$ endowed with
the split symplectic structure.
\end{itemize}
Let $\G$ be a finitely generated group
containing a {\bf U}-element $x$. Then
$\phi(x) = \done$ for every
homomorphism $\phi:\G \to G$.
\end{corollary}

\medskip
\noindent

For $\G = B(q,p)$ with $p {\big |} q$ and $G=
\Ham(M,\ome)$ this result can be extended
to all closed manifolds with $\pi_2 = 0$.

\begin{proposition}\label{prop1.6F}
Assume that $p$ divides $q$.
Let $(M,\ome)$  be an arbitrary closed
symplectic manifold with $\pi_2(M)=0$.  Then
$\phi(a)=\done$  for every homomorphism
$\phi:BS(q,p)\to \Ham (M,\ome)$.
\end{proposition}

\noindent
This result is easier than the previous ones.
It is proved in 2.6 below.

\medskip
\noindent
{\bf 1.6.I. Proof of Corollary 1.1.D:}

Let $G$ be the group of all symplectic diffeomorphisms
of a closed oriented surface $M$, and let $\G$ be
a non-uniform lattice as in 1.6.F.
Consider any homomorphism
$\phi : \G \to G$. It is proved in \cite{FM} (cf. \cite{KM})
that there exists a normal subgroup $K \subset \G$ of finite index
such that $\phi (K) \subset \Symp_0 (M)$.
Denote by $x$ a {\bf U}-element of $\G$ (see \cite{LMR} and 1.6.F).
Choose $p \in \NN$ such that
$x^p \in K$. It follows from  Proposition 2.2 of \cite{LMR} that
$x^p$ is a {\bf U}-element in $K$. Applying Corollary 1.6.G we conclude
that $\phi(x^p) = \done$. Therefore the kernel $Q$ of $\phi$
contains an element of infinite order,
and thus $Q \subset \G$ is an infinite normal subgroup. By Margulis
finitness theorem we get that $Q$ is of finite index in $\G$,
and hence $\phi$ has finite image.
\hfill $\square$

\medskip
\noindent
{\bf 1.6.J.} Let $\G$ be a non-uniform irreducible lattice as
in 1.6.F. Let $(M,\omega)$ be a closed symplectic manifold
with $\pi_2 = 0$. Suppose that the symplectic filling function
$v(s)$ of $M$ satisfies $v(s) \geq c s^{\epsilon}$ for some
$c,\epsilon > 0$. (Think, for instance, about the standard
symplectic torus $\TT^{2m}$).
Write $G = \Symp_0 (M,\omega)$ and $G_0 = \Ham(M,\omega)$.
We claim that {\it every homomorphism $\phi:\G \to G$ has finite image.}
Here is the proof.
It is known (see e.g. \cite{Ba}, \cite{MS}, \cite{LMP}) that there exists
a countable subgroup $E \subset H^1(M,\RR)$ and a homomorphism
\footnote{The subgroup $E$ is called the flux subgroup
of $(M,\omega)$. It is known \cite{LMP} to be discrete provided $\pi_2(M) = 0$,
though we do not use it.
With the notation of Subsection 2.2 below $E= \text{Image}(\Delta)$.
The homomorphism $\overline{\Flux}$ is defined in 2.2 in the simplest
case when $E = \{0\}$. In the general case one modifies it in the obvious way.}
$$\overline{\Flux}: G \to H^1(M,\RR)/E$$
whose kernel equals $G_0$.
It follows from Margulis finitness theorem that
the kernel $K$ of $\overline{\Flux}\circ \phi$
is a normal subgroup of finite index in $\G$.
Note that $\phi (K)$ is contained in $G_0$.

Denote by $x$ a {\bf U}-element of $\G$ (see \cite{LMR} and 1.6.F).
Choose $p \in \NN$ such that
$x^p \in K$. It follows from  Proposition 2.2 of \cite{LMR} that
$x^p$ is a {\bf U}-element in $K$. Applying 1.6.D we conclude
that $\phi(x^p) = \done$. Therefore the kernel $Q$ of $\phi$
contains an element of infinite order,
and thus $Q \subset \G$ is an infinite normal subgroup. By Margulis
finitness theorem we get that $Q$ is of finite index in $\G$,
and hence $\phi$ has finite image. The claim follows.
\hfill $\square$

\medskip
\noindent
{\bf 1.6.K.} The previous claim is in general not true
when one replaces the group $\Symp_0(M,\omega)$ with $\Diff_0(M)$.
One can give a counterexample already when $M$ is the 2-torus $\TT^2$.
Consider the group $\G = PSL(2,\ZZ[\sqrt{2}])$ of projectivized
matrices with determinant $1$ and with the entries of the form
$a+b\sqrt{2}$ where $a,b \in \ZZ$.
Call a number $a-b\sqrt{2}$ to be conjugate to $a+b\sqrt{2}$
in $\ZZ[\sqrt{2}]$. Given a matrix $A \in \G$, denote by $\bar A$ the
matrix with conjugate entries. Consider a monomorphism
$$\psi : \G \to PSL(2,\RR) \times PSL(2,\RR)\;, \;\;A \mapsto (A,{\bar A})\;.$$
One can show (see 2.12 in \cite{LMR} and \cite{vdG}) that $\psi(\G)$ is an irreducible
non-uniform lattice
in the Lie group $PSL(2,\RR) \times PSL(2,\RR)$. The real rank of this
Lie group equals 2. Therefore it follows from 1.6.J that every
homomorphism $\G \to \Symp_0(\TT^2)$ has finite image. On the other hand,
$\G$ embeds to $\Diff_0(\TT^2)$ in the obvious way:
$$\G \to PSL(2,\RR)\times PSL(2,\RR) \to \Diff_0(S^1) \times \Diff_0(S^1) \to
\Diff_0(S^1 \times S^1)\;.$$

\medskip
\noindent
{\bf 1.6.L.} Let us emphasize that the group
of Hamiltonian diffeomorphisms of a closed
symplectic manifold with $\pi_2(M)=0$ has no torsion.
This follows immediately from Theorem 1.6.A
(see also the paragraph following Proposition 2.6.A
below for a direct proof).
If in addition $\pi_1(M)$ has trivial center
then the same is valid for the group $\Symp_0(M,\omega)$
(use Theorem 1.6.B).
Looking at groups $\Ham(S^2)$ and $\Symp_0(\TT^2)$
which contain torsion elements we conclude that the topological
assumptions cannot be removed.

\section{A review of the symplectic action}\label{sec:review}

In this section we sum up some known facts on
the symplectic action which will be used for the
proof of results stated in the introduction.  Unless
otherwise stated, all symplectic manifolds
below are assumed to be connected.

\subsection{Action difference}\label{subsec:action}
Let $(P,\Ome)$  be a symplectic manifold with
$\pi_1(P)=\pi_2(P)=0$  (and hence $P$  is
necessarily non-closed).  Let $\varphi:P\to P$
be a symplectic diffeomorphism.  Given two fixed
points $x$  and $y$  of $\varphi$, define
their {\it action difference} $\del (\varphi; x,y)$
 as follows. Take any curve $\gam:[0;1]\to P$
with $\gam(0)=x$, $\gam (1)=y$  and take a disc
$\Sig\subset P$  with $\part\Sig =\varphi\gam
-\gam$  (here $\gam$ is considered as a 1-chain
in $P$).  Put
\renewcommand{\theequation}{\thesubsection.\Alph{equation}}
\setcounter{equation}{0}
\begin{equation}\label{eq2.1A}
\del (\varphi;x,y)=\intl_\Sig\Ome\ .
\end{equation}
Let us verify that this definition is correct,
that is $\del (\varphi ;x,y)$  does not depend
on the choice of $\gam$ and $\Sig$.  Indeed, let
$\gam ',\Sig'$  be another choice.  Since $P$ is
simply connected, there exists a disc $\Del$
with $\part\Del =\gam'-\gam$.  Note that the 2-chain
$\Pi =\Sig -\Sig' +\varphi\Del -\Del$
represents a 2-sphere in $P$, and hence $\int_\Pi\Ome
=0$  since $\pi_2(P)=0$.  But this yields
$$\intl_\Sig\Ome -\intl_{\Sig'}\Ome
=\intl_\Del\Ome -\intl_{\varphi\Del}\Ome =0$$
since $\varphi$  preserves $\Ome$, and the claim
follows.  The action difference behaves nicely
under iterations of $\varphi$.

\setcounter{theorem}{1}
\begin{proposition}\label{prop2.1B}
$\del (\varphi^n; x,y)=n\del(\varphi;x,y)$ for
all $n\in\ZZ$.
\end{proposition}

\pr Assume for simplicity that $n>0$.  Take a
curve $\gam$  joining $x$  with $y$, and let
$\Sig$  be a disc with $\part\Sig =\varphi\gam
-\gam$.  Put $\Del =\Sig +\varphi\Sig+\cdots
+\varphi^{n-1}\Sig$.  Clearly $\Del$  is a disc
with $\part\Del =\varphi^n\gam -\gam$.  Then
$\del (\varphi^n;x,y)=\intl_\Del\Ome
=n\intl_\Sig\Ome =n\del(\varphi; x,y)$.\hfill$\square$

\medskip
\noindent
The main result of the present section is as
follows.

\begin{theorem}\label{theo2.1C} Let $(M,\ome)$ be a
closed symplectic manifold with $\pi_2(M)=0$,
and let $f\in \Symp_0(M,\ome)$  be a symplectic
diffeomorphism with a fixed point of
contractible type.  Then $f$  admits a lift
$\tilf$ to the universal cover $(\tilM,\tilome)$
of $(M,\ome)$  such that $\del(\tilf;x,y)\not=
0$ for some fixed points $x$  and $y$  of $\tilf$.
\end{theorem}

\medskip
\noindent
The proof is given in 2.5 below.

\begin{remark}\label{remark2.1D}
If $f$ as above is  Hamiltonian, the lift $\tilf$
coincides with the canonical lift of $f$
defined in 1.5.
\end{remark}

\subsection{Symplectic and Hamiltonian}
\label{subsec:hamiltonian}
Let $(M,\ome)$ be a symplectic manifold (not
necessarily closed) with $\pi_2(M)=0$. Consider
a path $\{ f_t\}_{t\in [0;1]}$  of symplectic
diffeomorphisms with $f_0=\done$, $f_1=f$.  Let
$\xi_t$  be the corresponding time-dependent
vector field on $M$:
$$\frac{d}{dt}f_tx =\xi_t(f_tx)\ \text{for
all}\ x\in M\ ,\ t\in [0;1]\ .$$
Since the Lie derivative $L_{\xi_t}\ome$
vanishes we get that the 1-forms
$\lam_t=-i_{\xi_t}\ome$  are closed. Write
$[\lam_t]$  for the cohomology class of $\lam_t$.
The quantity
\setcounter{equation}{0}
\begin{equation}\label{eq2.2A}
\Flux \{ f_t\} =\intl^1_0[\lam_t]dt\in
H^1(M,\RR)
\end{equation}
is called the flux of the path $\{ f_t\}$.  A
path $\{ f_t\}$  is called {\it Hamiltonian\/}
if the 1-forms $\lam_t$  are exact for all $t$.
In this case there exists a smooth function
$F:M\times [0;1]\to\RR$ so that $\lam_t=dF_t$,
where $F_t(x)$  stands for $F(x,t)$. The
function $F$  is called the Hamiltonian function
generating the flow $\{f_t\}$. Note that $F_t$
is defined uniquely up to an additive
time-dependent constant.

A symplectic diffeomorphism $f:M\to M$  is
called Hamiltonian if there exists a Hamiltonian
path $\{ f_t\}_{t\in [0;1]}$ with $f_0=\done$
and $f_1=f$.  Hamiltonian diffeomorphisms form a
group denoted by $\Ham (M,\ome)$.
The next statement is proved in \cite{Ba},\cite{MS}.

\setcounter{theorem}{1}
\begin{proposition}\label{prop2.2B}
Let $(M,\ome)$  be a closed symplectic manifold.
Let $\{ f_t\}_{t\in [0;1]}$ be a path of
symplectic diffeomorphisms with $f_0=\done$  and
$\Flux \{ f_t\}=0$.  Then the diffeomorphism
$f_1$  is Hamiltonian.
\end{proposition}

\medskip
\noindent
It is well known (see \cite{Ba},\cite{MS}) that
$\Flux \{ f_t\}$  does not change under a homotopy of the
path $\{ f_t\}$  with fixed end points.  Thus
one can define a homomorphism
$$\Del :\pi_1(\Symp_0(M,\ome))\to H^1(M,\RR)$$
by $\Del(a)= \Flux \{ f_t\}$, where $\{ f_t\}$  is
a loop $(f_0=f_1=\done)$  of symplectic
diffeomorphisms representing an element
$a \in \pi_1(\Symp_0(M,\ome))$.  Sometimes $\Del$
vanishes identically. In this case, consider a map
\setcounter{equation}{2}
\begin{equation}\label{eq2.2C}
\overline{\Flux}: \Symp_0(M,\ome)\to H^1(M,\RR)\ ,
\end{equation}
which sends a diffeomorphism $f\in \Symp_0(M,\ome)$
to $\Flux \{ f_t\}$, where $\{ f_t\}$  is any
symplectic path with $f_0=\done$, $f_1=f$.
The condition $\Del\equiv 0$ guarantees that
$\overline{\Flux}$  is well defined.  Moreover,
it is a group homomorphism.  We refer to
\cite{MS} for a detailed discussion of the flux.
For the proof of Theorem \ref{theo1.6B} we shall need
the next result.

\setcounter{theorem}{3}
\begin{proposition}\label{prop2.2D}
Assume that $(M,\ome)$  is a
closed symplectic manifold with $\pi_2(M)=0$.
Suppose in addition that the fundamental group
$\pi_1(M)$ has trivial center.  Then $\Del$
vanishes and hence the homomorphism
$$\overline{\Flux}:\Symp_0 (M,\ome)\to H^1(M,\RR)$$
is well defined.
\end{proposition}

\medskip
\noindent
{\bf Proof:}
This fact can be easily extracted e.g. from \cite{LMP}.
For the reader's convenience we present the argument.
Take any loop $\{f_t\}$ of symplectic
diffeomorphisms representing an element
$a \in \pi_1\big(\Symp_0(M,\ome)\big)$.
Fix a point, say $y_*$, on $M$ and denote by $b$ the orbit
$\{f_ty_*\}$. Let $\beta \in \pi_1(M,y_*)$ be the
element represented by $b$.

 Let $c:[0;1] \to M$ be any closed curve
with $c(0)=c(1)= y_*$. Write $\gamma$
for the element represented by $c$
 in $\pi_1(M,y_*)$ and $\bar \gamma$ for the
homology class of $c$ in $H_1(M,\ZZ)$.

Consider the map
$[0;1]\times[0;1] \to M$,
$$(s,t) \mapsto f_t(c(s)).$$
It defines a 2-torus, say $\Sigma$, in $M$.
Since $b$ and $c$ lie on the 2-torus the elements
$\beta$ and $\gamma$ commute in $\pi_1(M,y_*)$.
This remains true for any choice of $c$. Therefore
$\beta$ belongs to the center of the fundamental group.
Since the center is trivial by our assumption, we conclude
that $\beta = 1$. In other words, the orbits of $\{f_t\}$
are contractible in $M$.

It is well known (see \cite{LMP},\cite{MS})
that
$$\langle \Del(a),{\bar\gamma} \rangle = - \intl_{\Sigma} \ome.$$
Filling the cycle $b \subset \Sigma$ by a 2-disc
we get that the torus $\Sigma$ is homologous
to a 2-sphere in $M$. Thus  the integral above vanishes
since $\pi_2(M)=0$. This proves that $\Del$ vanishes.
\hfill $\square$

\begin{remark}\label{remark2.2E}
It follows from \ref{prop2.2B} that
$\text{Ker} (\overline{\Flux})=\Ham (M,\ome)$.
\end{remark}

\begin{remark} The triviality of the center of $\pi_1$
has another consequence which will be used below.
Observe that exactly the same argument as in the proof
of 2.2.D shows that all orbits of any loop
of diffeomorphisms of $M$ are contractible in $M$.
Take any diffeomorphism $g\in \Diff_0
(M)$ and any path $\{g_t\}$  with $g_0=\done, g_1=g$.
Consider the lift $\{ \tilg_t\}$  of this path to the
universal cover $\tilM$  so that $\tilg_0=\done$.
The observation above implies that the
lift $\tilg_1$  does not depend on the choice of the
path $\{g_t\}$.  In particular, every $g\in \Diff_0(M)$
has a canonical lift $\tilg$  to $\tilM$. Of course, if
$g$ is a Hamiltonian diffeomorphism then $\tilg$ coincides
with its canonical lift defined in 1.5 above.
\end{remark}

\subsection{Symplectic action}\label{subsec:sympact}
Let $(M,\ome)$  be a symplectic manifold with
$\pi_2(M)=0$.  Let $\{f_t\}$  be a Hamiltonian path with
$f_1=f$  generated by a Hamiltonian function $F:M\times
[0;1]\to\RR$.
Let $x$  be a fixed point of $f$ such that its orbit
$\alp =\{f_tx\}_{t\in [0;1]}$ is contractible in $M$.
Take any 2-disc $\Sig\subset M$  with $\part\Sig =\alp$,
and define the symplectic action
\setcounter{equation}{0}
\begin{equation}\label{eq2.3A}
\calA(F, x)=\intl_\Sig\ome -\intl^1_0F_t(f_t x)dt\ .
\end{equation}
Since  $\pi_2(M)=0$ the integral
$\intl_\Sig\ome$  does not depend on the choice of the
disc $\Sig$.
The following deep fact is proved in \cite{Schwarz} by using
Floer homology.

\setcounter{theorem}{1}
\begin{proposition}\label{prop2.3B}
Let $(M,\ome)$  be a closed symplectic manifold with
$\pi_2(M)=0$. Let $\{f_t\}$, $f_0=\done$, $f_1=f$ be a
Hamiltonian path on $M$  generated by a Hamiltonian
function $F$.  Assume that $f\not=\done$.  Then $f$  has
a pair of fixed points $x$ and $y$ so that their orbits
$\{f_tx\}$  and $\{f_ty\}$  are contractible and
$\calA(F,y)-\calA(F,x)\not= 0$.
\end{proposition}

\medskip
\noindent
This proposition is the key ingredient from ``hard"
symplectic topology we use in this paper.

\subsection{Action difference revisited}
\label{subsec:revisited}

Let us return to the situation described in
\ref{subsec:action} above.  Let $(P,\Ome)$  be a
symplectic manifold with $\pi_1(P)=\pi_2(P)=0$.  Note
that any path $\{\varphi_t\}$ of symplectic
diffeomorphisms of $P$  is automatically Hamiltonian
since $H^1(P,\RR)=0$.  Take such a path and write $\Phi$
for the Hamiltonian function.  Write
$\varphi=\varphi_1$.

\begin{proposition}\label{prop2.4A}
$\del(\varphi;x,y)=\calA(\Phi,y)-\calA(\Phi,x)$.
\end{proposition}

\medskip
\noindent
This justifies the wording ``action difference".

\pr Consider the orbits $\alp_x=\{\varphi_tx\}_{t\in
[0;1]}$ and $\alp_y=\{\varphi_ty\}_{t\in [0;1]}$  of $x$
and $y$, and choose discs $\Sig_x,\Sig_y$  in $P$  so
that $\part\Sig_x=\alp_x$  and $\part\Sig_y=\alp_y$.
Choose a curve $\gam:[0;1]\to P$ with
$\gam(0)=x$  and $\gam(1)=y$.  Define a
2-chain $\Del:[0;1]\times [0;1]\to P$ by
$\Del(t,s)=\varphi_t\gam(s)$.  Note that $\part\Del =-\gam
+\varphi\gam -\alp_y+\alp_x$, where we assume that the
boundary of the square $[0;1]\times [0;1]$  is oriented
counter-clockwise.  Thus the boundary of the topological
disc $\Pi =\Del +\Sig_y-\Sig_x$  equals $\varphi\gam
-\gam$.  Therefore
$$\del (\varphi ;x,y)=\intl_\Pi\Ome\ .$$
Denote by $\xi_t$  the vector field of the flow
$\varphi_t$  (see \ref{subsec:hamiltonian} above).  Then
\begin{eqnarray*}
&&\Del^*\Ome =\Ome\Big(
\xi_t(\varphi_t\gam(s)),\frac{\part}{\part
s}\varphi_t\gam (s)\Big)dt\wedge ds\\[0.3em]
&&\qquad =-d\Phi_t\Big(\frac{\part}{\part
s}\varphi_t\gam(s)\Big)dt\wedge ds\ .
\end{eqnarray*}
Hence
\begin{eqnarray*}
&&\intl_\Del\Ome =\intl_{[0;1]\times [0;1]}\Del^*\Ome
=-\intl^1_0dt\intl^1_0ds\ d\Phi_t\Big(\frac{\part}{\part
s}\varphi_t\gam (s)\Big)\\[0.3em]
&&\qquad =\intl^1_0\Phi_t(\varphi_tx)dt-\intl^1_0\Phi_t
(\varphi_ty)dt\ .
\end{eqnarray*}
Therefore,
$$\del(\varphi;x,y)=\intl_\Pi\Ome =\intl_\Del\Ome
+\intl_{\Sig_y}\Ome -\intl_{\Sig_x}\Ome
=\calA(\Phi,y)-\calA(\Phi,x)\ .$$
The proof is complete.\hfill $\square$

\subsection{Proof of Theorem 2.1.C}
\label{subsec:proofoftheorem}
The proof splits into two cases.

\medskip
\noindent
{\bf Case I: $f$ is Hamiltonian.}
Let $\{ f_t\}$  be a Hamiltonian path on $M$  with
$f_0=\done$, $f_1=f$.  Denote by $F$  the Hamiltonian
function.  Let $\tilf_t,\tilf,\tilF$ be the lifts of
$f_t,f$,$F$  to the universal cover $(\tilM,\tilome)$
respectively.

Proposition \ref{prop2.3B} guarantees that $f$  has two
fixed points $x,y$  such that their orbits are contractible
and $\calA(F,y)-\calA(F,x)\not= 0$. Let $\tilx,\tily$  be
any lifts of $x$ and $y$  to $\tilM$.  The
contractibility of the orbits yields $\tilf\tilx =\tilx$  and
$\tilf\tily =\tily$.  Further,
$$\calA(\tilF,\tilx)=\calA(F,x)\quad {\rm
and}\quad\calA(\tilF,\tily)=\calA(F,y)\ .$$
Combining this with \ref{prop2.4A} we get
$$\del(\tilf,x,y)=\calA(\tilF,y)-\calA(\tilF,x)=\calA(F,y)
-\calA(F,x)\not= 0\ .$$
This proves the statement of the theorem in Case I.

\medskip
\noindent
{\bf Case II: $f$ is not Hamiltonian.}  Let $x\in M$
be a fixed point of contractible type of $f$.  Thus
there exists a path of symplectic diffeomorphisms $\{
f_t\}$  such that $f_0=\done$, $f_t=f$  and the orbit
$\{ f_tx\}$  is contractible.  Note that $\Flux \{ f_t\}
\not= 0$ in view of Proposition \ref{prop2.2B}.  Thus
there exists an element $\alp\in\pi_1(M,x)$  such that
$\langle \Flux \{ f_t\}, \oalp\rangle\not= 0$, where
$\oalp$  stands for the image of $\alp$  in $H_1(M,\ZZ)$
under the Hurewitz homomorphism.

Let $\{\tilf_t\}$  be the lift of $\{f_t\}$  to the
universal cover $\tilM$. Put $\tilf = \tilf_1$.
Denote by $\tilF$  the
Hamiltonian function of $\{\tilf_t\}$  (recall from
\ref{subsec:revisited} that every path of symplectic
diffeomorphisms on $\tilM$ is Hamiltonian).  Note that
$d\tilF_t=\tau^*\lam_t$  where $\tau:\tilM\to M$  stands
for the natural projection, and $\{\lam_t\}$ is the
family of 1-forms associated to $\{f_t\}$  as in
\ref{subsec:hamiltonian} above.  Let $T:\tilM\to\tilM$
be the deck transformation corresponding to
$\alp\in\pi_1(M,x)$.  We claim that
\setcounter{equation}{0}
\begin{equation}\label{eq2.5A}
\tilF_t(Tz)-\tilF_t(z)=\langle
[\lam_t],\oalp\rangle\quad\text{for all}\quad t\in
[0;1],z\in\tilM\ .
\end{equation}
Indeed, choose any path $\gam:[0;1]\to\tilM$ with
$\gam(0)=z$  and $\gam(1)=Tz$.  Then
$$\tilF_t(Tz)-\tilF_t(z)=\intl_\gam d\tilF_t=\intl_\gam
\tau^*\lam_t=\intl_{\tau(\gam)}\lam_t=\langle
[\lam_t],\oalp\rangle\ ,$$
and \eqref{eq2.5A} follows.

Let $\tilx$  be a lift of $x$.  Then $\tilf\tilx =\tilx$
since $\{ f_tx\}$  is
contractible.
Moreover, $\tilf_t$
commutes with $T$  for every $t$  so that $\tilf T\tilx
=T\tilx$.  We claim that
\begin{equation}\label{eq2.5B}
\del (\tilf ;\tilx,T\tilx)=-\langle \Flux \{
f_t\},\oalp\rangle\;.
\end{equation}
Indeed,
$$\del (\tilf,\tilx,T\tilx)=\calA(\tilF,T\tilx)
-\calA(\tilF,\tilx)\ .$$
Choose a disc $\Sig\subset\tilM$ spanning
$\{\tilf_t\tilx\}$.  Then $T\Sig$ spans
$\{\tilf_tT\tilx\}$.  Applying
\eqref{eq2.3A} we get that
\begin{eqnarray*}
&&\calA(\tilF,T\tilx)-\calA(\tilF,\tilx)=
\intl_{T\Sig}\tilome -\intl_\Sig\tilome\\[0.3em]
&&\qquad -\intl^1_0\Big[\tilF_t(T\tilf_t\tilx)
-\tilF_t(\tilf_t\tilx)\Big] dt\\[0.3em]
&&\overset{2.5.A}{=\!=\!=}-\intl^1_0\langle
[\lam_t],\oalp\rangle dt\overset{2.2.A}{=\!=\!=}
-\langle \Flux \{ f_t\},\oalp\rangle\ .
\end{eqnarray*}
This proves \eqref{eq2.5B}.  Recall now that $\langle
\Flux \{ f_t\},\oalp\rangle \not= 0$.  Hence $\del
(\tilf;\tilx,T\tilx)\not= 0$  as required. This completes
the proof of Theorem
\ref{theo2.1C}.\hfill$\square$

\subsection{Action spectrum of Hamiltonian
diffeomorphisms}\label{subsec:spectrum}
Let $(M,\ome)$  be a closed symplectic manifold with
$\pi_2(M)=0$.  For a Hamiltonian diffeomorphism $f\in
\Ham (M,\ome)$  take a Hamiltonian path $\{f_t\}$  of
symplectic diffeomorphisms with $f_0=\done$, $f_1=f$.
Let $F(x,t)$  be the corresponding Hamiltonian function
normalized so that
$$\int_M F(x,t) d(\text{volume})=0$$ for all
$t\in [0;1]$.  Let $\text{Fix}_0$  be the set of all fixed
points $x$ of $f$ such that the orbit $\{f_tx\}$  is
contractible.  Define the action spectrum of $f$  as
$$\spectrum (f)=\{\calA (F,x)\mid x\in
\text{Fix}_0\}\subset\RR\ .$$
It is known (see \cite{Schwarz}) that the action
spectrum of $f$ does not depend on the choice of a
Hamiltonian path generating $f$, and is invariant under
conjugations in the group of all symplectomorphisms
of $M$.
Moreover, this set is
compact \cite{HZ},\cite{Schwarz}.  Define the following
invariant:
$$\width (f)=\underset{\alp,\bet\in \spectrum (f)}
{\max\mid\alp -\bet|}\ .$$
Lifting the path $\{f_t\}$ to the universal cover $\tilM$
and applying Proposition 2.4.A we get
$$\width (f) = \max_{x,y} \delta(\tilf,x,y),$$
where $\tilf$ is the canonical lift of $f$ to $\tilM$
and $x,y$ run over the set of fixed points of $\tilf$.
Combining this with \ref{prop2.1B}, 2.1.C and 2.1.D  we get the
following inequality.

\begin{proposition}\label{prop2.6A}
$\width (f^n)\succeq n$ for every $f\in \Ham (M,\ome)\bks
\{\done\}$.
\end{proposition}

\medskip

This fact turns out to be very useful for studying
group-theoretic properties of $\Ham(M,\omega)$.
For instance, let us show that
for a closed symplectic manifold $(M,\omega)$ with $\pi_2=0$
the group
$\Ham(M,\omega)$ has
no torsion. Indeed, suppose that $f$ is a Hamiltonian diffeomorphism
of finite order. Then the sequence $\width (f^n)$ is bounded, and so
Proposition 2.6.A yields $f=\done$.

As another immediate application, we prove {\rm Proposition} {\rm
\ref{prop1.6F}}.

\medskip
\noindent
{\bf Proof of 1.6.H:}
Assume first that
that $f^m$  is
conjugate to $f$  in $\Ham(M,\ome)$  for some $m\in\ZZ$
with $|m|>1$. Then $f^{m^k}$ is conjugate to $f$ for
all $k \in \NN$. Hence
$$\width(f^{m^k}) = \width(f)$$
for all $k \in \NN$ which contradicts Proposition 2.6.A.
Thus we proved the proposition for $B(m,1)$.

Consider now the general case of $B(q,p)$ where
$p$ divides $q$. Assume that $q=pm$
with $|m|>1$
and
$f^{pm}$ is conjugate to $f^p$.
Applying the argument above to $f^p$
we see that $f^p = \done$. The abscence
of torsion yields $f=\done$ as required.
\hfill $\square$

\section{Proofs of lower bounds for
growth}\label{sec:bounds}

In 3.2--3.4  below we prove Theorems \ref{theo1.1B},
\ref{theo1.4A} and \ref{theo1.5B}.

\subsection{Remarks on the symplectic filling function}
\label{subsec:function}
We use notations of subsection \ref{subsec:symplectic}.

\begin{lemma}\label{lemma3.1A}
Symplectic filling functions $v_1$  and $v_2$
associated to Riemannian metrics $\rho_1$  and $\rho_2$
on $M$  and to the same base point $x\in\tilM$  are
equivalent:  $v_1\sim v_2$.
\end{lemma}

\pr There exists $c>1$ such that
$$B_{\tilrho_2}(c^{-1}s)\subset B_{\tilrho_1}(s)\subset
B_{\tilrho_2}(cs)$$
and
$$c^{-1}|\xi|_{\tilrho_2}\le |\xi|_{\tilrho_1}\le
c|\xi|_{\tilrho_2}$$
for all $s\ge 0$, $\xi\in T^*\tilM$.
Put $u_i(s)=\inf\limits_{\alp\in\calL}\
\sup\limits_{z\in B_{\tilrho_i}(s)}|\alp_z|_{\rho_i}$
and $w_i(s)=su_i(s)$ for $i=1,2$.  Then
$c^{-1}u_1(c^{-1}s)\le u_2(s)\le cu_1(cs)$  for all $s$,
and so
$$w_1(c^{-1}s)\le w_2(s)\le w_1(cs)\ .$$
Since $v_1,v_2$  are inverse to $w_1(s),w_2(s)$
respectively we conclude that $c^{-1}v_1(s)\le v_2(s)\le
cv_1(s)$  as required.\hfill $\square$

\begin{lemma}\label{lemma3.1B}
Symplectic filling functions $v_1$ and $v_2$  associated
to Riemannian metrics $\rho_1$  and $\rho_2$ and to
base points $x_1,x_2\in\tilM$  are equivalent.
\end{lemma}

\pr The group $\Symp_0(M,\ome)$  acts transitively on
$M$. Hence there exists a symplectomorphism
$\varphi:\tilM\to\tilM$  which commutes with the action
of $\pi_1(M)$  on $\tilM$  so that $\varphi(x_1)=x_2$.
Consider the symplectic filling function $v_3(s)$
associated to the metric $\rho_3=\varphi^*\rho_1$  and to
the base point $x_2$.  Note that $\varphi$  establishes
a 1-to-1 map $\calL\to\calL$  which sends $\alp$ to
$\varphi^*\alp$.  Hence $v_3(s)\equiv v_1(s)$.  But
$v_3\sim v_2$  in view of \ref{lemma3.1A}.  This
completes the proof.\hfill $\square$

\begin{lemma}\label{lemma3.1C}
$v(cs)\sim v(s)$  for all $c>0$.
\end{lemma}

\pr It suffices to show that
$$v(s)\le v(cs)\le cv(s)$$
for all $s>0$, $c>1$.
The inequality $v(cs)\ge v(s)$  holds since $v$  is
increasing.  By definition, $v(s)u(v(s))=s$. Hence
$$cv(s)=\frac{cs}{u(v(s))}$$
and
$$v(cs) = \frac{cs}{u(v(cs))}\ .$$
But $v(cs)\ge v(s)$, so $u(v(cs))\ge u(v(s))$ since $u$
is non-decreasing. It follows that
$$v(cs)\le\frac{cs}{u(v(s))}=cv(s)\ $$
as claimed.
\hfill $\square$

\subsection{Starting the proofs}
\label{subsec:proofs}
Let $(M,\ome)$  be a closed symplectic manifold with
$\pi_2(M)=0$.  Let $f\in \Symp_0(M,\ome)\bks\{\done\}$  be
a symplectic diffeomorphism with a fixed point of
contractible type.  Choose its lift $\tilf$  as in
Theorem \ref{theo2.1C}.  Let $x,y\in\tilM$ be fixed
points of $\tilf$  with non-vanishing action difference:
$|\del (\tilf;x,y)|=c>0$.  Then \ref{prop2.1B} yields
\setcounter{equation}{0}
\begin{equation}\label{eq3.2A}
|\del(\tilf^n;x,y)|=nc\quad\text{for all}\quad
n\in\NN\ .
\end{equation}
Take any curve $\gam:[0;1]\to\tilM$ joining $x$  and $y$
on $\tilM$.  Denote by $\ell_n$ the loop formed by
$\gam$  and $\tilf^n\gam$.  Then

\begin{equation}\label{eq3.2B}
nc= {\Big |}\intl_{\ell_n}\alp \;{\Big |}
\end{equation}
for every primitive $\alp$ of $\tilome$.
 Denote by $b$ the length of $\gam$. Clearly,
\begin{equation}\label{eq3.2C}
\text{length} (\ell_n)\le b(1+\Gam_n(f))\le 2b\Gam_n(f)
\end{equation}
(by definition, $\Gam_n(f)\ge 1$).
Without loss of generality assume that the
base point appearing in the definition of the
symplectic
filling function $v$ is
$x$.

\subsection{Proof of 1.4.A and 1.5.B}\label{subsec:threethree}

Denote by $B(s)$ the ball of
radius $s$  centered at $x$.  Pick any positive $s>b$.

\medskip
\noindent
{\bf Case I.}  Assume that $\tilf^n(\gam)\subset
B(s)$.  Then \ref{eq3.2B} yields
$$nc\le \text{length} (\ell_n) \cdot
\sup\limits_{z\in B(s)}|\alp_z|\;.$$
Taking into account \ref{eq3.2C} together with the
definition of $u$  we get $nc\le 2b\Gam_n(f)u(s)$.  Thus

\setcounter{equation}{0}
\begin{equation}\label{eq3.3A}
\Gam_n(f)\ge\frac{nc}{2bu(s)}\;.
\end{equation}

\medskip
\noindent
{\bf Case II.} Assume
that there exists a point $z\in\gam$  so that $\tilf^n(z)$
lies outside $B(s)$.  Then obviously
\begin{equation}\label{eq3.3B}
\Gam_n(f)\ge\frac{s}{b}\ .
\end{equation}
Given $n$, the choice of $s$ is in our hands.  Choose it
(assuming that $n$  is large enough) so that
$$\frac{s}{b}=\frac{nc}{2bu(s)},$$ that is
$s=v\Big(\frac{1}{2}nc\Big)$.  Then \eqref{eq3.3A} and
\eqref{eq3.3B} yield
$$\Gam_n(f)\ge\frac{s}{b}=\frac{1}{b}v\Big(\frac{1}{2}
nc\Big)\ .$$
Applying Lemma \ref{lemma3.1C} we get $\Gam_n(f)\succeq
v(n)$. This completes the proof of Theorem 1.4.A.

\medskip
\noindent
Assume now that $f$  is a
Hamiltonian diffeomorphism which does not propagate.
Choose $s_0>0$  large enough and argue as above.
The only possible case {\it for all\/}
$n\in\NN$  is Case I.  Thus
\eqref{eq3.3A} reads
$$\Gamma_n(f)\ge\frac{nc}{2bu(s_0)}\succeq
n,$$
which proves Theorem 1.5.B.
\hfill $\square$

\subsection{Proof of 1.1.B
(following \cite{PS})}

Here we assume that $M=\TT^2$.  An additional ingredient
is an isoperimetric inequality which was proved (in a
much more general context) by Bonk and Eremenko
\cite{BE}.  Let $\RR^2$  be the Euclidean plane and
$L\subset\RR^2$  a lattice.  There exists a constant
$\kap =\kap (L)>0$  such that for every
piece-wise smooth
curve $\bet:
S^1\to\RR^2\bks L$ which is contractible in
$\RR^2\bks L$
\setcounter{equation}{0}
\begin{equation}\label{eq3.4A}
\Area (\bet)\le\kap\cdot \Length (\bet)\ .
\end{equation}
Here
$\Area(\bet)=\inf_\varphi\int_{D^2}|\varphi^*\tilome |$,
where $\tilome$  is the Euclidean area form and
$\varphi$  runs over all piece-wise smooth
maps $D^2\to\RR^2\bks L$  with
$\varphi\Big|_{S^1}=\bet$.  We refer to \cite{PS} for a
different proof of \eqref{eq3.4A}.
In order to prove \ref{theo1.1B} start arguing as in 3.2.
Further, choose a vector $e\in\ZZ^2$ and
a natural number $N$  so that the curves $\gam$  and
$\tilf\gam$ are homotopic with fixed end points in
$\RR^2\bks L$,  where the lattice $L$  is defined by
$L=x+e+N\ZZ^2$. Since $L$ consists of fixed points of
$\tilf$  we see that $\tilf^n\gam$  is homotopic to
$\tilf^{n-1}\gam$ with fixed end points in $\RR^2\bks
L$.  Therefore the loop $\ell_n$  is contractible in
$\RR^2\bks L$  which yields
\begin{equation}\label{eq3.4B}
\Area (\ell_n)\le \kap \cdot \Length (\ell_n)\ .
\end{equation}
Combining this with \ref{eq3.2B} and \ref{eq3.2C}
we get
$$nc={\Big|}\intl_{\ell_n}\alp {\Big |}\le \Area (\ell_n)\le\kap\cdot
\Length (\ell_n)\le 2\kap b\Gam_n(f)\ .$$
We conclude that $\Gam_n(f)\succeq n$.  This completes
the proof.\hfill $\square$

\section{Proofs of lower bounds for distortion}
\label{sec:slow}

\subsection{Measurements on the group of symplectomorphisms}
\label{subsec:meas}

Let $(M,\omega)$ be a closed symplectic manifold.
Throughout \S 4 we fix a compatible Riemannian metric
on $(M,\omega)$, that is a metric of the form
$\omega(\xi,J\eta)$ where $J$ is an almost complex structure on $M$.
An important feature of such a metric is the equality
$|\nabla F | = |\xi|$, where $F$ is any smooth function on $M$
with the Hamiltonian field $\xi$. Let $\{f_t\},\; t\in [0;1]$
be a smooth path of symplectomorphisms. Write $\xi_t$ for the time-dependent
vector field on $M$ which generates this path.
Put
$$L\{f_t\} = \int_0^1 \max_{x \in M} |\xi_t(x)| dt.$$
If $\{f_t\}$ is a Hamiltonian path, write $F(x,t)=F_t(x)$ for
its Hamiltonian function. We always assume that the Hamiltonian
function is normalized so that the mean value of $F_t$ with respect to the
canonical measure on $M$ vanishes for all $t \in [0;1]$.
Put
$$\Lambda\{f_t\}= \int_0^1 \max_{x \in M} |F_t(x)| dt.$$

Both $L$ and $\Lambda$ are right-invariant length structures
on groups $\Symp_0 (M,\omega)$ and $\Ham (M,\omega)$ respectively
associated to norms $\max |\xi|$ and $\max |F|$ on their
Lie algebras. The right-invariance, of course, means
that
$L\{f_th\} = L\{f_t\}$ and
$\Lambda\{f_th\} = \Lambda\{f_t\}$.
In fact $\Lambda$ is left-invariant as well and is a version
of Hofer's length structure on the group of Hamiltonian
diffeomorphisms (see e.g. \cite{Pbook}). The length structure $L$
has an ${\bf L_2}$-cousin which is used in hydrodynamics. Let us
emphasize that, as it should be for length structures, both $L$
and $\Lambda$ do not change under reparameterization of paths,
and are additive with respect to juxtaposition.

Given $f\in \Symp_0(M,\omega)$ set
$$\alpha (f) = \inf L\{f_t\},$$
where the infimum is taken over all symplectic paths $\{f_t\}$
with $f_0=\done$ and $f_1 = f$.

Given $f \in \Ham (M,\omega)$ set
$$\beta (f) = \inf {\big (} L\{f_t\} + \Lambda\{f_t\}{\big )},$$
where the infimum is taken over all Hamiltonian paths $\{f_t\}$
with $f_0 = \done$ and $f_1 = f$.

One readily checks that both $$\alpha: \Symp_0(M,\omega) \to [0;+\infty)$$
and $$\beta: \Ham (M,\omega) \to [0;+\infty)$$
are {\it norms}. In other words, they have the
following properties:
\begin{itemize}
\item{}
triangle inequality:
$$\alpha(fg) \leq \alpha(f)+\alpha(g) \;\;\text{and}\;\;
\beta(fg) \leq \beta(f)+\beta(g);$$
\item{} symmetry:
$\alpha(f^{-1})=\alpha(f)$ and $\beta(f^{-1}) = \beta(f)$;
\item{} non-degeneracy: $\alpha(f) > 0$ and $\beta(f) > 0$
for all $f\neq \done$.
\end{itemize}
The non-degeneracy follows from the obvious inequalities
$$\alpha (f) \geq \max_{x \in M} \text{dist}(x,fx),$$
and
$$\beta (f) \geq \max_{x \in M} \text{dist}(x,fx).$$
Let us mention that the triangle inequality is the only
property of $\alpha$ and $\beta$ we need for the proof of
Theorems 1.6.A,B.

\subsection{Geometric inequalities}

Assume now that $\pi_2(M,\omega) = 0$. Let $u$ be the
function defined in subsection 1.2 and let $d$ be the
diameter of $M$.

\begin{lemma}
Let $f \in \Ham (M,\omega)$ be a Hamiltonian diffeomorphism
with $\beta (f) = b$. Then
$$\width (f) \leq 2(b+bu(d+b)).$$
\end{lemma}

\begin{lemma} Assume in addition that $\pi_1(M)$ has trivial
center. Let $f \in \Ham (M,\omega)$ be a Hamiltonian diffeomorhism
with $\alpha (f) = a$. Then
$$\width (f) \leq 2(a^2+da+au(d+a)).$$
\end{lemma}

\noindent
The proofs are postponed till 4.5.

\subsection{Proof of Theorem 1.6.A}

Let $\G \subset \Ham (M,\omega)$ be a finitely
generated group with generators $h_1,...,h_N$.
Put $\kap =\max_j \beta(h_j)$. Then $\beta (f) \leq \kap ||f||$
for all $f \in \G$. Put $c_n = ||f^n||$ and $w_n = \width (f^n)$.
Applying Lemma 4.2.A to $f^n$ we get that
$$w_n \leq 2\kap c_n +2\kap c_n u(d+\kap c_n).$$
Thus the exists a constant $\mu>0$ independent on $n$
so that
$$w_n \leq \mu c_n u(\mu c_n).$$
This yields
$$\mu c_n \geq v(w_n),$$
where $v$ is the symplectic filling function. Recall from
2.6.A that $w_n \succeq n$. Applying Lemma 3.1.C we conclude
that $c_n \succeq v(n)$. This completes the proof.
\hfill $\square$

\subsection{Proof of Theorem 1.6.B}

Let $\G \subset \Symp_0 (M,\omega)$ be a finitely
generated group with generators $h_1,...,h_N$.
Take $f \in \G$.

\medskip
\noindent
{\bf Case 1: $f$ is Hamiltonian.}
Put $\kap =\max_j \alpha(h_j)$. Then $\alpha (f) \leq \kap ||f||$
for all $f \in \G$. Put $c_n = ||f^n||$ and $w_n = \width (f^n)$.
Applying Lemma 4.2.B to $f^n$ we get that
$$w_n \leq 2\kap^2 c_n^2 +2d\kap c_n +2 \kap c_n u(d+\kap c_n).$$
Thus the exists a constant $\mu>0$ independent on $n$
so that
$$2w_n \leq \mu c_n^2 +\mu c_n u(\mu c_n).$$
Therefore for every $n$ either $\mu c_n^2 \geq w_n$ or
$\mu c_n  u(\mu c_n) \geq w_n$. The last inequality
is equivalent to $\mu c_n \geq v(w_n)$,
where $v$ is the symplectic filling function. Recall from
2.6.A that $w_n \succeq n$. Applying Lemma 3.1.C we conclude
that $$c_n \succeq \min {\big (}{\sqrt n}, v(n){\big )}.$$
This completes the proof for a Hamiltonian $f$.

\medskip
\noindent
{\bf Case 2: $f$ is not Hamiltonian.}  We shall
use the flux homomorphism $$\overline{\Flux}:\Symp_0(M,\ome)
\to H^1(M,\RR)$$ defined in \ref{eq2.2C},
\ref{prop2.2D}.  Assume that
$f^n=h_{i_m}\cdot\ldots\cdot h_{i_1}$.  Then
$$n \cdot \overline{\Flux}
(f)=\sum^m_{j=1} \overline{\Flux} (h_{i_j}).$$  Choose any norm on
$H^1(M,\RR)$  and put $$\kap =\max\limits_{j\in \{1;\ldots
;N\}}\mid \overline{ \Flux} (h_j)|.$$  We get
$$n|\overline{\Flux} (f)|\le m\kap,$$  so $$\| f^n\|\ge\frac{|
\overline{\Flux}
(f)|}{\kap}\cdot n.$$  Since $\overline{\Flux} (f) \not= 0$  we conclude
that $$\| f^n\|\succeq n.$$
This completes the proof.
\hfill $\square$

\subsection{Proof of geometric inequalities}

We start with the following situation which is common
for 4.2.A and 4.2.B.
Let $x_*\in\tilM$
be the base point of $\tilM$. Let $D\subset \tilM$ be the ball of radius $d$
centered at $x_*$. Note that $D$ projects onto the whole $M$.
Let $\{f_t\},\; t\in [0;1]$ be any path of symplectic diffeomorphisms
of $M$ with $f_0 = \done$ and $f_1 = f$. Lift it to the path $\{h_t\}$
on $\tilM$, and write $\xi_t$ for the vector field generating $h_t$.
Denote $h=h_1$.

Suppose that $h$ is the canonical Hamiltonian lift of $f$
defined in 1.5.
Recall that the path $\{h_t\}$ on $\tilM$ is always Hamiltonian
since $\tilM$ is simply connected.
Let $H(x,t)$ be any Hamiltonian function on $\tilM$ generating
$\{h_t\}$ (the choice of a time-dependent additive constant
is in our hands).
 It follows from Proposition 2.4.A that
$$ \width (f) = \max \delta(h;y_1,y_2) = \max (\calA(H,y_1) -\calA(H,y_2)),$$
where $y_1,y_2$ run over fixed points of $h$ and
$$\calA (H,y) = \int_{\{h_ty\}} \lambda - \int_0^1 H(h_ty,t) dt,$$
where $\lambda$ is a primitive of $\tilde \omega$.
Thus our purpose is to estimate
$\calA (H,y)$, where $y$ is a fixed point of $h$.

Note also that
$$\calA (H,y) - \calA (H,Ty) = \delta(h;y,Ty) = 0$$
for every deck transformation $T$ of the covering $\tilM \to M$.
To see this, consider
any Hamiltonian path on $M$ joining
$\done$ with $f$. Applying Proposition 2.4.A to the lift of
this path we obtain $\delta(h;y,Ty) = 0$.
Summing up this discussion, we can assume without loss
of generality that $y \in D$.

The ``symplectic area" term  of $\calA (H,y)$ can be estimated as
follows.
Denote by $\gam$
the orbit $\{h_ty\},\; t \in [0;1]$. Abbreviate $c = L\{f_t\}$.
Obviously
\setcounter{equation}{0}
\begin{equation}\label{4.5.A}
\Length (\gamma) = \int_0^1 |\xi_t(h_ty)| dt \leq c.
\end{equation}
Therefore $\gamma$ is contained in the ball of radius $d+c$
with the center at $x_*$.
Hence
\begin{equation}\label{4.5.B}
{\Big |}\int_{\gamma} \lambda {\Big |} \leq cu(d+c)
\end{equation}
for every primitive $\lambda$ of $\tilde \omega$.

It remains to estimate the ``Hamiltonian" term of $\calA(H,y)$.
At this point the proofs of 4.2.A and 4.2.B split up.

\medskip
\noindent
{\bf Proof of Lemma 4.2.A:} Suppose that $\{f_t\}$ is a Hamiltonian
path which joins $\done$ with $f$.
Therofore $h$ is the canonical lift of $f$.
Write $F(x,t)$ for the normalized
Hamiltonian of $\{f_t\}$, and choose $H$ to be
the lift of $F$ to $\tilM$.
Note that
$${\Big |}\int_0^1 H(h_tx,t) dt {\Big |} \leq \Lambda\{f_t\}.$$
Applying \eqref{4.5.B} we obtain
$$\calA (H,y) \leq \Lambda \{f_t\} + cu(d+c)$$
with $c = L\{f_t\}$.  Since this is true for every Hamiltonian
path $\{f_t\}$ which joins $\done$ with $f$ we get
$$\calA(H,y) \leq b+bu(d+b).$$
The same action bound holds for any other fixed point of $h$.
Therefore
$$\width (f) \leq 2(b+bu(d+b)),$$
which proves 4.2.A.
\hfill $\square$

\medskip
\noindent
{\bf Proof of Lemma 4.2.B:}
Here $\{f_t\}$ is an arbitrary symplectic path joining $\done$
with $f$.
The triviality of the center
of $\pi_1(M)$ guarantees that $h$ is the canonical lift
of $f$ (see Remark 2.2.F above). Suppose that the Hamiltonian
$H(x,t)$ is normalized so that $H_t(x_*) = 0$ for all $t\in [0;1]$.
{\bf Warning:} since the path $\{f_t\}$ is not in general Hamiltonian,
the function $H$ does not  come as
a lift of any function on $M$, and in particular
$H$ may be unbounded on $\tilM$.
We will get around this difficulty by noticing that
the differential of $H_t$ is a bounded 1-form.
We proceed as follows.
Observe that
$$|H_t(h_ty)| \leq |H_t(y)| +  \Length (\gamma) \cdot \max_x |\nabla H_t(x)|,$$
and
$$|H_t(y)| \leq |H_t(x_*)| + d\cdot \max_x |\nabla H_t(x)|.$$
Taking into account that $H_t(x_*) = 0,\;\;|\nabla F_t|= |\xi_t|$
and using (4.5.A) we get
$$|H_t(h_tx)| \leq (d+c)\max_x |\xi_t(x)|.$$
Therefore
$${\Big |} \int_0^1 H_t(h_tx) dt{\Big |} \leq (d+c)\cdot \int_0^1 \max_x |\xi_t(x)|dt =
(d+c)c.$$
Thus
$$\calA (H,y) \leq c^2+cd+cu(d+c).$$
 Since this is true for every symplectic
path $\{f_t\}$ which joins $\done$ with $f$ we get
$$\calA (H,y) \leq a^2+ad+au(d+a).$$
The same action bound holds for any other fixed point of $h$.
Therefore
$$\width (f) \leq  2(a^2+ad+au(d+a)),$$
which proves 4.2.B.
\hfill $\square$

\subsection{A remark on smoothness}
Here we outline an extension of our results above to
$C^1$-smooth symplectic diffeomorphisms.
Geometric inequalities 4.2.A,B
remain valid for symplectomorphisms which are generated
by $C^1$-smooth vector fields.
Various topological facts about symplectic diffeomorphisms
which appeared above should extend to the $C^1$-case without
problems. One should simply use an appropriate approximation
by $C^{\infty}$-diffeomorphisms.

A more delicate argument is needed to show
that any Hamiltonian diffeomorphism
$f\neq \done$ has two fixed points with strictly positive
action difference (see crucial Proposition 2.3.B).
We claim that this holds true for $C^1$-smooth Hamiltonian
diffeomorphisms. The proof of the claim is based on a remarkable
``energy-capacity" inequality in symplectic topology. We need a
version from \cite{Schwarz}. Let $g$ be a $C^{\infty}$-smooth
Hamiltonian diffeomorphism of $M$ which displaces an open subset
$B \subset M$:
$$g(B) \cap B = \emptyset.$$
Then $g$ has two fixed points whose
action difference is at least $c(B)$,
where $c(B)$ is a strictly positive constant
which depends only on $B$.
Assume now that $f \neq \done$ is
a $C^1$-smooth Hamiltonian diffeomorphism.
Clearly,
$f$ displaces a small ball, say $B$. Choose a sequence $g_i$ of
$C^{\infty}$-smooth Hamiltonian diffeomorphisms which converges to $f$
in the $C^1$-sense. Then $g_iB \cap B = \emptyset$ for large $i$.
Therefore each $g_i$ has a pair of fixed points, say $x_i$ and $y_i$,
whose action difference is at least $c(B)$. By compactness, choose
a subsequence $i_k \to \infty$ such that $x_{i_k}$ and $y_{i_k}$
converge to fixed points $x$ and $y$ of $f$ respectively. Since
the action difference is continuous
with respect to the $C^1$-convergence
(use 2.1 above to see this), we conclude that the action difference
of $x$ and $y$ is at least $c(B)$. This proves the claim.

\subsection{Proof of Theorem 1.1.E}
The proof is divided into 3 steps.

\medskip
\noindent
1)Let $G$ be the group of all symplectic diffeomorphisms
of a closed oriented surface $M$ of genus $\geq 2$.
Let $\phi: BS(q,p) \to G$ be a homomorphism. We assume
for simplicity that $q>p>0$. Denote by $\text{Mod}$
the mapping class group of $M$, and let $\pi:G \to \text{Mod}$
be the natural projection. Farb-Lubotzky-Minsky theorem
\cite{FLM} implies that every non-torsion element in $\text{Mod}$
is undistorted (in the sense of 1.6.C above). Since the element
$a \in BS(q,p)$ has logarithmic distortion (see 1.6.E) we conclude
that $\pi (\phi (a))$ is of finite order in $\text{Mod}$.
Therefore there exists $k \in \NN$ such that
$\phi (a^k)$ lies in $\Symp_0 (M)$. Denote
$f = \phi (b), g = \phi (a^k)$.
Then $g \in \Symp_0$ and
$g^q = fg^pf^{-1}$.

\medskip
\noindent
2) We claim that in fact $g$ lies in $\Ham (M,\omega)$.
Indeed, assume on the contrary that
$$\overline{\Flux} (g) \neq 0.$$
Consider the isomorphism of $H^1(M,\RR)$ induced by $f$,
and denote by $I$ its inverse. Then the equation
$g^q = fg^pf^{-1}$ yields
$$\overline{\Flux} (g^q) = I\cdot \overline{\Flux} (g^p).$$
Rewrite this as
$$I\cdot \overline{\Flux} (g) =
\frac{q}{p} \cdot \overline{\Flux} (g).$$
Therefore the number $q/p$ is a root of the characteristic
polynomial $\chi(t)$ of $I$. Note now that since $I$ preserves the
lattice $H^1(M,\ZZ)$ all coefficients of $\chi$ are integers.
Moreover the leading and the free coefficients of $\chi$ are equal
to $1$. Hence the only rational roots of $\chi$ are $\pm 1$,
so $q = \pm p$.
This contradicts to the assumption $q > p >0$.
The claim follows.

\medskip
\noindent
3) Note that
$$g^{q^m} = f^mg^{p^m}f^{-m}$$
for all $m \in \NN$.
Then
$$\width (g^{q^m}) = \width (g^{p^m}).$$
Assume that $g\neq\done$. Then 2.6.A yields
$$
\width (g^{q^m}) \succeq q^m.$$
Applying  4.2.A we get
$$\width (g^{p^m}) \preceq \beta(g^{p^m}) \preceq p^m.$$
Combining  these inequalities we conclude that
$q^m \preceq p^m$ which contradicts our choice of $p$ and $q$.
Therefore $g=\done$, and so $\phi(a)^k = \done$.
This completes the proof.
\hfill $\square$

\subsection{Epilogue }
It is instructive to take a look at our approach
to distortion from the viewpoint of global geometry
of symplectomorphisms groups.

Let $(M,\omega)$ be a closed symplectic manifold
with $\pi_2 = 0$ and such that $\pi_1$ has trivial
center. Choose a compatible Riemannian metric on $M$
and choose any norm on $H^1(M,\RR)$.
Define a norm $$\gamma:\Symp_0(M,\omega) \to [0;+\infty)$$
by
$$\gamma(f) = \alpha (f) + {\big |} \overline{\Flux} (f){\big |}.$$
Inequality 4.2.B in conjunction with 2.6.A gives
a bound for distortion of
any non-trivial cyclic subgroup of $\Symp_0(M,\omega)$ with respect
to $\gamma$.
Assume for instance that $(M,\omega)$ is symplectically hyperbolic.
Then $\gamma (f^n) \succeq n$ for every non-Hamiltonian diffeomorphism
$f$ and $\gamma (f^n) \geq {\sqrt n}$ for every Hamiltonian $f$.

Similarly, inequality 4.2.A leads to distortion bounds for
cyclic subgroups of $\Ham (M,\omega)$ with respect to the
norm $\beta$.

Of course, these distortion bounds give rise to obstructions
for reperesentations of finitely generated groups into
symplectomorphisms groups. Let me mention also that the group
of Hamiltonian diffeomorphisms carries in addition two
remarkable norms which are invariant under conjugations --
the Hofer norm and the commutator norm. It would be interesting
to understand distortion bounds in these cases. We refer
to \cite{Pbook} for some results on distortion in the Hofer norm.
As far as the commutator norm is concerned, we refer to
\cite{BaGh} and \cite{E}.

\section{Appendix: an example}

In this section we show that the growth bound 1.4.A is
in general not true for symplectic diffeomorphisms
which have fixed points but none of them
is of contractible type. Namely, we will construct
a diffeomorphism $f \in \Symp_0(M,\ome)$ of certain
symplectic manifold
$(M,\ome)$ with $\pi_2(M)=0$ which has the following
properties:
\begin{itemize}
\item [$\ast$] $f$ has fixed points, and none of them
is of contractible type;
\item [$\ast$] $f^2=\done$, so that $\{\Gam_n(f)\}$  is
a bounded sequence.
\end{itemize}

The construction goes as follows.
Consider the standard symplectic torus $\TT^4
=\RR^4/\ZZ^4$  endowed with the symplectic form
$dp_1\wedge dq_1+dp_2\wedge dq_2$. Let
$$\gam:\TT^4\to\TT^4,\;\;
\gam(p_1,q_1,p_2,q_2)=(p_1,q_1+\frac{1}{2},-p_2,-q_2)\;\mod\;1$$
be an involution of the torus. Clearly, $\gam$ is
symplectic and has no fixed points.  Thus it generates a
free symplectic action of $\ZZ_2$  on $\TT^4$. Let
$(M,\ome)$  be the quotient space $\TT^4/\ZZ_2$, and let
$\tau:\TT^4\to M$  be the natural projection.

Consider the flow
$$\of_t:\TT^4\to\TT^4\;,\;\;(p_1,q_1,p_2,q_2)\to
(p_1,q_1+\frac{t}{2},p_2,q_2)\;\mod\;1.$$
It is a symplectic flow
which commutes with $\gam$.  Hence it defines a
symplectic flow $f_t$  on $M$. Put $f=f_1$  and
note that $f^2 =\done$. Look now at
fixed points of $f$. All of them have the form
$\tau(z)$,  where $z\in\TT^4$  satisfies $\of_1z =\gam
z$, so that
$$z=(p,q,m_1/2,m_2/2)\;\mod\;1\;,\; \;p,q\in\RR,\;
m_1,m_2\in\ZZ.$$
  We claim that all fixed points of $f$
are {\it not\/} of contractible type.  To prove this,
assume on the contrary that there exists another path,
say $\{ g_t\}_{t\in [0;1]}$  of symplectic
diffeomorphisms of $M$  such that $g_0=\done$, $g_1=f$
and the loop $\{ g_t\tau(z)\}$  is contractible.
Consider a loop $\{ h_t\}$ of diffeomorphisms of $M$
formed by $\{ f_t\}$  and $\{ g_{1-t}\}$:
$h_tx=f_{2t}x$ if $t\in [0;1/2]$  and $h_t x=g_{2-2t}x$
if $t\in [\frac{1}{2};1]$.  Let $\{ \oh_t\}$  be its
lift to $\TT^4$.  Note that $\tau\circ\oh_1=\tau$,  so
either $\oh_1=\done$  or $\oh_1=\gam$.  Due to our
construction one has $\oh_1z =\gam z$, which rules out
the first possibility.  Therefore $\oh_1=\gam$.  But
$\gam$  acts non-trivially on $H_1(\TT^4)$  and thus
cannot be isotopic to the identity. This contradiction
proves the claim.

\vfill\eject
\noindent
{\bf Acknowledgments.} My interest to finitely
generated groups of symplectomorphisms originated
from a number of conversations with Marc Burger and Matthias Schwarz.
I am greatly indebted to Marc Burger for his comments
on the preliminary version of this paper:
He pointed out the explicit link to
the Zimmer
program and in particular
explained to me Corollary 1.1.D and
its reduction to Theorem 1.1.C.
I thank Dima Burago, Alex Eremenko,
Anatole Katok, Shahar Mozes, Zlil Sela,
Misha Sodin and Zeev Rudnick
for useful discussions, and Jean-Claude
Sikorav
for his permission to include a result of preprint \cite{PS}
into the present paper. I am grateful to Felix Schlenk
for careful reading the manuscript and many useful critical
remarks.


\begin{thebibliography}{AAAAA}

\bibitem[DAG]{DAG}
G.~D'Ambra and M.~Gromov,
Lectures on transformations groups: geometry and dynamics,
in {Surveys in Differential Geometry}, supplement to the
{\it J. Diff. Geom.}, number 1, 1991, pp. 19-112.

\bibitem[Ba]{Ba} A.~Banyaga, Sur la structure du groupe des
diff\'eomorphismes qui pr\'eservent une forme symplectique,
{\it Comm. Math. Helv.}
{\bf 53} (1978), 174-227.

\bibitem[BG]{BaGh}J.~Barge and E.~Ghys, Cocyles d'Euler et de
Maslov, {\it Math. Ann.} {\bf 294}(292), 235-265.

\bibitem[BPS]{BPS} P.~Biran, L.~Polterovich and D.~Salamon,
Propagation in Hamiltonian dynamics and relative symplectic
homology, preprint math.SG/0108134.

\bibitem[BE]{BE} M.~Bonk and A.~Eremenko,
Uniformly hyperbolic surfaces, {\it Indiana University
Math. Journal}, {\bf 49} (2000), 61-80.


\bibitem[BM]{BurMon} M.~Burger and N.~Monod, Bounded cohomology
of lattices in higher rank Lie groups, {\it J. Eur. Math. Soc.}
{\bf 1} (1999), 199-235.

\bibitem[E]{E} M.~Entov, Commutator length for symplectomorphisms,
preprint math.SG/0112012.

\bibitem[FF]{FF} B.~Farb and J.~Franks,
Groups of homeomorphisms of one-manifolds I: actions of non-linear
groups, preprint math.DS/0107085.

\bibitem[FLM]{FLM} B.~Farb, A.~Lubotzky and Y.~Minsky, Rank-1
phenomena for mapping class groups, {\it Duke Math. J.} {\bf 106}
(2001), 581-597.

\bibitem[FM]{FM} B.~Farb and H.~Masur,
Superrigidity and mapping class groups,
{\it Topology} {\bf 37}(1998), 1169-1176.

\bibitem[FS]{FS} B.~Farb and P.~Shalen, Real-analytic actions of lattices,
{\it Invent. Math.} {\bf 135} (1999), 273-296.


\bibitem[Fl]{Fl} A.~Floer, Morse theory for Lagrangian intersection
theory, {\it J. Diff. Geom.} {\bf 28} (1988), 513-517.

\bibitem[G1]{GJDG} M.~Gromov, K\"ahler hyperbolicity and $L^2$-Hodge
theory, {\it J. Diff. Geom.} {\bf 33} (1991), 263-292.

\bibitem[G2]{Gromov} M.~Gromov, Asymptotic invariants of infinite
groups, in {\it Geometric group theory,} Vol.2 (Sussex 1991), 1-295.
LMS Lecture Note Ser., 182, Cambridge Univ. Press, 1993.

\bibitem[Gh]{Ghys} E.~Ghys, Actions de reseaux sur le cercle,
{\it Invent. Math.} {\bf 137} (1999), 199-231.

\bibitem[HaK]{HK} B.~Hasselblatt and A.~Katok,
Principle structures, Preprint, to appear in {\it Handbook
of Dynamical Systems}.

\bibitem[HZ]{HZ} H.~Hofer and E.~Zehnder, {\it Symplectic invariants and
Hamiltonian dynamics}, Birkh\"auser Advanced Texts,
Birkh\"auser Verlag, 1994.


\bibitem[KM]{KM} V.~Kaimanovich and H.~Masur, The Poisson boundary
of the mapping class group, {\it Invent. Math.}
{\bf 125} (1996), 221-264.

\bibitem[LMP]{LMP} F.~Lalonde, D.~McDuff
and  L.~Polterovich, On the flux
conjectures,
in {\it ''Geometry, Topology and Dynamics, ed F. Lalonde,
Proceedings of the CRM 1995 Workshop
in Montreal, the CRM Special Series of the AMS''},
Vol.15, 1998, pp. 69-85.

\bibitem[LMR]{LMR} A.~Lubotzky, S.~Mozes and M.S.~Raghunathan,
The word and riemannian metrics on lattices in semisimple Lie
groups, IHES Publ. Math. {\bf 91} (2000), 5-53.

\bibitem[MS]{MS} D.McDuff and D.Salamon, {\it Introduction
to symplectic
topology}, Oxford Mathematical Monographs,
Oxford University Press, 1995.

\bibitem[P1]{Pbook} L.~Polterovich,
{\it The geometry of the group of
 symplectic diffeomorphisms}, Lectures in Mathematics ETH Z\"urich,
 Birkh\"auser, 2001.

\bibitem[P2]{P2} L.~Polterovich, Dissipation in contact
dynamics, to appear in {\it Ergod. Theory and Dynam. Syst.},
revised version of preprint math.SG/0009227.

\bibitem[PSi]{PS} L.~Polterovich and J.-C.~Sikorav, A
linear isoperimetric inequality for the punctured
Euclidean plane, preprint math.GR/0106216.

\bibitem[PSo]{PSo} L.~Polterovich and M.~Sodin,
Diffeomorphisms with an anomalous growth of the
differential, preprint math.DS/0110306.

\bibitem[Si]{Sikorav}  J.-C.~Sikorav, Growth of a primitive
of a differential
form, {\it Bull. Soc. Math. France}, to appear.

\bibitem[Sch]{Schwarz} M.~Schwarz, On the action spectrum
for closed symplectically aspherical manifolds,
{\it Pacific J. Math.} 193(2000), 419-461.

\bibitem[vdG]{vdG} G.~van der Geer, {\it Hilbert modular surfaces},
Springer, 1988.

\bibitem[Z]{Zimmer} R.~Zimmer, Actions of semisimple groups and discrete
subgroups, {\it Proc. I.C.M.}, Berkeley 1986, p.1247-1258.


\end{thebibliography}
\end{document}